\documentclass[12pt]{article}

\usepackage{amssymb}
\usepackage{amsmath}
\usepackage{cleveref}

\usepackage{subcaption}
\newcommand{\mm}[1]{{\boldsymbol{{ #1}} }}
\newcommand{\vv}[1]{{\boldsymbol{{ #1}} }}

\usepackage{newtxtext,newtxmath}
\usepackage{graphicx}

\usepackage[letterpaper,margin=1in]{geometry}

\renewenvironment{abstract}
	{\quotation}
	{\endquotation}

\date{}

\makeatletter
\renewcommand{\fnum@figure}{\textbf{Figure \thefigure}}
\renewcommand{\fnum@table}{\textbf{Table \thetable}}
\makeatother

\usepackage{url}

\def\scititle{
	Approximating Analytic Spectra of Hyperbolic Systems with Summation-by-Parts Finite Difference Operators}
\title{\bfseries \boldmath \scititle}

\author{
	Brittany A. Erickson$^{1}$\and
	\small$^{1}$Department of Computer Science, University of Oregon, Eugene Oregon, 97402, USA.
}

\begin{document} 

% Insert the title and author list
\maketitle

% Abstract, in bold
% There are strict length limits, and not all formats have abstracts.
% Consult the journal instructions to authors for details.
% Do not cite any references in the abstract.
\begin{abstract} \bfseries \boldmath
% Start with one or two sentences of background
In this work we explore the fidelity of numerical approximations to the analytic spectra of hyperbolic partial differential equation systems with variable coefficients. We are particularly interested in the ability of discrete methods to accurately discover sources of physical instabilities. By considering the perturbed equations that arise in linearized problems, we study systems in which a lower-order term can act as a source of internal energy within the system. We apply high-order accurate summation-by-parts finite difference operators, with weak enforcement of boundary conditions through the simultaneous-approximation-term technique, which leads to a provably stable numerical discretization with formal order of accuracy given by $p = 2, 3, 4$ and $5$. We derive analytic solutions using Laplace transform methods, which provide important ground truth to ensure numerical convergence at the correct theoretical rate. We derive the analytic spectrum and find that it is better captured with mesh refinement, although dissipative strict stability (where the growth rate of the discrete problem is bounded above by the analytic) is not obtained. We also find that sole reliance on mesh refinement can be a problematic means for determining physical growth rates as some eigenvalues emerge (and persist with mesh refinement) based on spatial order of accuracy but are non-physical. We suggest that numerical methods be used to approximate the spectra when numerical stability is guaranteed and convergence of the numerical spectra is evident with both mesh refinement and increasing order of accuracy.

\end{abstract}

\section{Introduction}
Systems of hyperbolic partial differential equations emerge across vast application areas, from gravitational waves to drug delivery targeting cancerous tissues, to unsteady volcanic eruptions \cite{PhysRevD.58.044020, FERREIRA201881, Karlstrom2016}.  
Determining the physical growth and/or decay rates involved in such processes is the topic of many studies, for example \cite{Zapadka2017, Blanco31122022, HALE2007163, Yarushina2015}, as they reveal physical mechanisms for unstable phenomena.  Such instabilities threaten, for example, the structural integrity of buildings and airplanes, and may indicate unstable dynamical growth that leads to regime changes in fluid systems.  

When nonlinearities are present, the governing equations can be notoriously challenging to analyze and solve and are often linearized around a background and/or steady state solution. Even so, analytic solutions can be difficult to obtain and one often turns to numerical methods in order to explore the underlying physics and/or origin of unstable regimes. However, there are several stumbling blocks when using numerical methods to explore physical instabilities.  If one observes an instability in the numerical simulation it can be difficult to determine if it is numerical or physical in origin. In addition, the spectrum of the discretization could appear to converge to eigenvalues that are artifacts of the numerical method but are non-physical.

In this work we explore the fidelity of numerical approximations to analytic spectra, focusing on variable-coefficient hyperbolic systems containing a lower-order term, which arise in perturbation studies (e.g. when nonlinear equations are linearized, for example in \cite{Karlstrom2016}). Using the Laplace transform, we derive the analytic spectra for several parameter regimes, as well as analytic solutions which provide ground truth for our numerical studies. We utilize a class of high-order-accurate finite difference methods satisfying a summation-by-parts (SBP) property \cite{kreiss_finite_1974,kreiss_existence_1977,strand_summation_1994}. SBP methods together with a weak enforcement of boundary conditions through the simultaneous-approximation-term (SAT) technique provide a framework for provably stable numerical discretizations, see for example \cite{DELREYFERNANDEZ2014171,SVARD201417,RANOCHA201820}, which will provide the assurance that instabilities are not numerical in origin.

% The outline is not required, but we show an example here.
The paper is organized as follows. Related work is discussed in \Cref{sec:related}. The governing equations considered are outlined in 
\Cref{sec:main}, with the analytic spectra derived in \Cref{sec:cont_spec}. Details of the SBP-SAT numerical discretization is given in  \Cref{sec:discrete}, with numerical verification and derivation of analytic solutions in \Cref{sec:verification} and numerical spectra in \Cref{sec:discrete_spec}. Conclusions follow in
\Cref{sec:conclusions}.

\section{Related Work}\label{sec:related}
Existence and uniqueness conditions for initial-boundary-value problems involving first-order hyperbolic systems with constant coefficients were derived in \cite{Hersh1963}, with analytic solutions constructed via Laplace transform methods. Well-posed boundary conditions for the corresponding variable coefficient problems were addressed in \cite{Kreiss1970} with results discussed in detail in \cite{Higdon1986}. More recently, \cite{Hastir2024} focused on a class of linear hyperbolic systems similar to those studied in this work (although limited to the special case where information travels solely in one direction) and well-posed boundary conditions and analytic solutions are derived.

Numerical studies include \cite{Hughes2014}, for example, where accuracy of finite-element approximations to the analytic spectrum were compared to NURBS (Non-Uniform Rational B-Splines) methods for elliptic, parabolic and hyperbolic problems. Studies focused on SBP-SAT discretizations have been explored in \cite{Nord2006, BergNordstrom2012, Erickson2019}, where first-order hyperbolic problems with variable coefficients were considered. In these works, the vector system was assumed to be solely coupled at the domain boundaries, and analytic and numerical spectra were calculated. Variable coefficients giving rise to strictly stable schemes were identified in \cite{Nord2006}, with physical growth explored in \cite{Erickson2019} and internal interface effects studied in \cite{BergNordstrom2012}.

In this work we extend these studies, exploring SBP-SAT approximations to analytic spectra associated with first-order, linear hyperbolic systems with variable coefficients, two-way information flow, and component coupling within the domain (through a lower order term) as well as at domain boundaries. In addition, we derive closed-form analytic solutions.

\section{Governing Equations}
\label{sec:main}
We consider hyperbolic systems of the form 
\begin{equation}\label{eqn: physical}
    U_t + A(x) U_x = B(x) U + F(x, t), \quad x \in [0, L], \quad t \geq 0,
\end{equation}
where $F$ is a known source function, $A(x), B(x) \in \mathbb{R}^{2 \times 2}$ are spatially varying coefficient matrices. Here $A(x)$ is diagonalizable with real eigenvalues and $B(x)U$ is a lower order term often referred to as the ``damping" term, when it does not contribute to growth of the system's internal energy \cite{Phillips1959}.  We assume $A = V\Lambda V^{-1}$, with $\Lambda = \text{diag}(-\bar{c}(x), \bar{c}(x))$, which pertains to applications with information flowing in two directions at speed $\bar{c}(x)$, for example the Euler equations in fluid dynamics where $\bar{c}(x)$ is the fluid sound speed. $U = [u_1, u_2]^T$ is the vector of physical variables. For the 1D Euler equations, for example, $U$ would correspond to the vector of fluid velocity and pressure. 

%In most applications, boundary conditions are desired on the physical variables: for example consider the following Dirichlet and Neumann boundary conditions (respectively)
%\begin{subequations}\label{eqn: physical_bc}
%\begin{align}
%    u_1(0, t) &= g_0(t),\\
%    \frac{\partial u_2}{\partial x}(L, t) &= g_L(t),
%    \end{align}
%\end{subequations}
%where $g_0, g_L$ are the known boundary data. 

An energy estimate for \eqref{eqn: physical} can be obtained through the diagonalization
\begin{equation}\label{eqn: char}
    W_t + \Lambda(x) W_x = \tilde{B}(x) W + \tilde{F}(x, t), \quad x \in [0, L], \quad t \geq 0,
\end{equation}
where $W = V^{-1}U$, $\tilde{B} = -\Lambda V^{-1} V_x + V^{-1}B V$, $\tilde{F} = V^{-1}F$. We consider the general case where $\tilde{B}(x) = [a(x), b(x); c(x), d(x)] \in \mathbb{R}^{2 \times 2}$ and derive conditions in which $\tilde{B}$ contributes to internal energy growth. Boundary conditions for \eqref{eqn: char} are motivated by the energy estimate we derive shortly and are given by
\begin{subequations}\label{eqn: char_bc}
\begin{align}
    w_2(0, t) &= R_0 w_1(0, t) + h_0(t),\\
   w_1(L, t) &= R_L w_2(L, t) + h_L(t),
    \end{align}
\end{subequations}
where $R_0, R_L$ are reflection coefficients constrained by the energy estimate.  Note that a non-diagonal $\tilde{B}$ matrix in \eqref{eqn: char} couples the components within the interior of the domain, whereas conditions \eqref{eqn: char_bc} couple components at boundaries for $R_0, R_L \neq 0$.

%For example, non-reflecting/absorbing boundary conditions correspond to $R_0 = R_L = 0$.  Dirichlet or Neumann conditions on the physical variables are imposed by setting reflection coefficients to either 1 or -1, for example, boundary conditions \eqref{eqn: physical_bc} are imposed by setting $R_0 = 1, R_L = -1$ and taking $h_0(t) = 2g_0(t)$, $h_L(t) = 2g_L(t)$. Other values for $R_0, R_L$ correspond to Robin boundary conditions.

Supplementing \eqref{eqn: char}-\eqref{eqn: char_bc} with the initial condition \begin{equation}\label{eqn: char_ic}
    W(x, 0) = f(x),
\end{equation} taking $\tilde{F} = h_0 = h_L = 0$ (not required but eases the analysis) an energy estimate can be obtained by multiplying \eqref{eqn: char} by $W^T$, adding the transpose, and integrating by parts. This process yields the estimate 
\begin{align}\label{eqn: energydot}
    \dot{E} &= \bar{c}(L)(R_L^2 - 1)w_2^2(L) + \bar{c}(0)(R_0^2 - 1)w_1^2(0) + \notag \\
    &(W, \Lambda_x W) + (W, \tilde{B}W) + (\tilde{B}W, W), 
\end{align}
where we define the energy $E = ||W||^2$ using the $L_2$ scalar product and norm
\begin{equation}
    (U, W) = \displaystyle\int_0^L U^T W dx, \quad ||U|| = \sqrt{(U, U)}.
\end{equation}
Estimate \eqref{eqn: energydot} illustrates the well-known result that reflection coefficients should satisfy $-1 \leq R_0, R_L \leq 1$ for there to be no energy growth from the boundaries. Assuming this and letting $\alpha = ||\tilde{B}|| = \underset{||U|| = 1}{\text{max}}||\tilde{B} U||$ and $\gamma = ||\Lambda_x||$,
applying the Cauchy-Schwarz inequality to \eqref{eqn: energydot} further reduces it to 
    \begin{align}\label{eqn: energydot2}
            \dot{E} &\leq (W, \Lambda_x W) +  (W, \tilde{B}W) + (\tilde{B}W, W) \notag\\
            & \leq (\gamma + 2\alpha) E.
    \end{align}
Equation \eqref{eqn: energydot2} can be integrated in time to yield the final result 
\begin{equation}
    E(t) \leq \exp\left[(\gamma + 2\alpha) t\right]E(0) = \exp\left[(\gamma + 2\alpha) t\right]||f||^2,
\end{equation}
which shows that energy in the system is bounded in terms of the data of the problem, i.e. the problem is well posed.  Some exponential growth must be tolerated, however, due to the effects of $\tilde{B}$ and $\Lambda_x$ (the latter of which includes those from spatially varying wave speeds), see chapter 3 in \cite{GKO}. 

\section{Analytic Spectra}
\label{sec:cont_spec}
Specific rates of growth or decay associated with the initial boundary value problem \eqref{eqn: char}-\eqref{eqn: char_ic} can be determined by computing the 
analytic spectrum.  We denote the Laplace transform of a locally integrable function $h(t)$ by 
\begin{equation}\label{eqn: laplace}
    \mathcal{L}[h] = \displaystyle\int_0^\infty h(t) \exp\left(-st\right) dt, \quad s \in \mathbb{C}.
\end{equation}
Taking the initial data $f = 0$ (and still considering $h_0, h_L, F = 0$ as in the energy estimate) and Laplace transforming in time yields
\begin{equation}\label{eqn: char1}
    s\hat{W} + \Lambda(x) \hat{W}_x = \tilde{B}(x) \hat{W}, \quad x \in [0, L], \quad s \in \mathbb{C},
\end{equation}
where we use a hat to denote the Laplace transform. Equation \eqref{eqn: char1} can be re-written as 
\begin{equation}\label{eqn: char2}
   \hat{W}_x = M(x, s)  \hat{W}, \quad x \in [0, L], \quad s \in \mathbb{C},
\end{equation}
where \begin{align}
M(x, s) &= \Lambda^{-1}(x)(\tilde{B}(x) - sI) \notag \\ &= \begin{pmatrix}-(a(x)-s)/\bar{c}(x)&-b(x)/\bar{c}(x)\\c(x)/\bar{c}(x)&(d(x)-s)/\bar{c}(x)\end{pmatrix}.
\end{align}
In this work we compute the associated analytic spectrum, considering several cases for the spatially varying coefficients.  Case 1 assumes constant coefficients and that $\tilde{B}$ is a diagonal matrix. This is a particular case of the more general, but is useful for illustrative purposes.  In Case 2 we still assume constant coefficients, but assume a  slightly more general case where $\tilde{B}$ is not diagonal. We then include some results for special cases of variable coefficients, namely Case 3 when $\tilde{B}$ is still a diagonal matrix.  For these three cases the solution to 
\eqref{eqn: char2} is given by 
\begin{equation}\label{eqn: general}
    \hat{W}(x, s) = \exp\left(\int^x M(\xi, s) d\xi\right)C(s),
\end{equation}
where $C(s) = [C_1(s), C_2(s)]$ are coefficients determined by the boundary conditions. We also include the set-up for Case 4 (the most general case) when $\tilde{B}$ is not diagonal, which will be explored in numerical studies. 

%The fundamental theorem of (linear autonomous) ODEs
%states that if $\lambda$ is an eigenvalue of multiplicity $m$ of $M$, then letting $v_1, v_2$ be $m$ linearly independent solutions of $(M - \lambda I)^m v = 0$ and $L(x, \lambda) = I + x(M - \lambda I)$, then $Y_k(x) = e^{\lambda x} L(x,\lambda)v_k, \quad k = 1, ..., m$ are $m$ linearly independent solutions of \eqref{eqn: char2} \cite{principles}. 

\subsection{Case 1: Constant coefficients; $\tilde{B}$ is a diagonal matrix}
In this case we have that $b = c = 0$, thus the eigenvalues of $M$ are given by \begin{equation}\lambda_{1, 2}(s) = -(a-s)/\bar{c}, \quad (d-s)/\bar{c},
\end{equation}
with eigenvectors $v_1, v_2$ the standard basis vectors of $\mathbb{R}^2$.  Thus \eqref{eqn: general} reduces to
\begin{equation}\label{eqn: sol}
  \hat{W}(x, s) = C_1(s)\exp\left[\lambda_1 (s)x\right]v_1 + C_2(s)\exp\left[\lambda_2(s) x\right]v_2.  
\end{equation} Substituting the Laplace transform of boundary conditions \eqref{eqn: char_bc} into \eqref{eqn: sol} yields the linear system $A(s)C(s) = 0$ where \begin{equation} A(s)=\begin{pmatrix}-R_0&1\\\exp\left[-(a-s)L/\bar{c}\right]&-R_L\exp\left[(d-s)L/\bar{c}\right]\end{pmatrix}.
\end{equation}
As we are interested in non-trivial solutions to \eqref{eqn: char2}, we must have that $\text{det}(A) = 0$, i.e. 
\begin{equation}\label{eqn: cond}
R_0 R_L \exp\left[L\lambda_2(s)\right]- \exp\left[L\lambda_1(s)\right] = 0.
\end{equation}
The above condition \eqref{eqn: cond} is met at discrete $s_n \in \mathbb{C}$ which define the analytic spectrum, namely 
\begin{equation}\label{eqn: s_n}
s_n = \frac{(a + d)L + \bar{c}\ln(R_0 R_L) + 2\bar{c}\pi i n}{2L}, \quad n \in \mathbb{Z}.
\end{equation}
Note that this result implies that the spectrum only exists if $R_0 R_L \neq 0$, namely, the two characteristic variables must be coupled at domain boundaries; a non-reflecting boundary condition will imply no analytic spectrum.  Furthermore, values of $s_n$ with positive real part can exist if $a+d$ is sufficiently positive. For example, if $R_0 = R_L \in \{ -1, 1 \}$, then $\text{Re}(s_n) > 0$ if $a+d>0$, leading to physical growth of the system. However, the term $\ln(R_0R_L)$ will counteract positive values of $a+d$ if it is sufficiently negative, leading to energy decay.

\subsection{Case 2: Constant coefficients; $\tilde{B}$ is not diagonal}
Here the eigenvalues of $M$ are 
\begin{equation}
    \lambda_{1, 2}(s) = \frac{-(a-d)}{2\bar{c}} \pm \frac{1}{2}\sqrt{\text{disc}(M)},
\end{equation}
where the discriminant of $M$ is \begin{equation}
    \text{disc}(M) = \left(\frac{a-d}{\bar{c}}\right)^2 - 4\left[\frac{bc - (a-s)(d-s)}{\bar{c}^2}\right].
\end{equation}
We first look for $s$-values corresponding to distinct eigenvalues (where $\text{disc}(M) \neq 0$, i.e. $s$-values where $s \neq \frac{a+d}{2} \pm \sqrt{bc})$, each with multiplicity 1. Since by assumption both $b$ and $c$ cannot both be zero, we assume WLOG that $b \neq 0$. Then the linearly independent eigenvectors of $M$ corresponding to each eigenvalue are given by 
\begin{equation}
    v_{1,2} = [1, \quad [-(a-s) - \lambda_{1,2}\bar{c}]/b].
\end{equation}
In this case, $\text{det}(A) = 0$ yields the equation
\begin{align}\label{eqn: detA}
    &\exp\left(\lambda_1 L\right)(-1 + R_L v_1^{(2)})(-R_0 + v_2^{(2)}) +\notag \\ &\exp\left(\lambda_2 L\right)(1-R_L v_2^{(2)})(-R_0 + v_1^{(2)}) = 0.
\end{align}
If $R_0 = R_L = R \in \{-1, \,\, 1\}$ then \eqref{eqn: detA} simplifies to 
\begin{equation}\label{eqn: root2}
    [\exp\left(\lambda_1 L\right) - \exp\left(\lambda_2 L\right)](-1 + R v_1^{(2)})(-R + v_2^{(2)})  = 0.
\end{equation}
We note that by inspection, the second two terms in \eqref{eqn: root2} yield the real root 
\begin{equation} 
s_0 = \frac{a + R(b+c) + d}{2}.
\end{equation}
Complex roots also exist when \begin{equation}\label{eqn: hard}
\exp\left[\lambda_1(s) L\right] = \exp\left[\lambda_2(s) L\right],
\end{equation}
(although by assumption, $\lambda_1 \neq \lambda_2$)
%If $\text{disc}(M) \in \mathbb{R}$ (although $\text{disc}(M) \neq 0$ as we assume distinct eigenvalues of $M$) we take the log of both sides of \eqref{eqn: hard} which yields $2|\text{disc}(M)| L = 2 \pi i n, \quad n \in \mathbb{Z}$, 
and one can check that solutions to \eqref{eqn: hard} are the roots
\begin{equation}
    s^\pm_n = \frac{a+d}{2} \pm \sqrt{bc - \pi^2 n^2 \bar{c}^2/L^2}, \quad n \in \mathbb{Z} \symbol{92} \{0\},
\end{equation}
which produces purely real roots for all $n \in \mathbb{Z}\symbol{92} \{0\}$ such that 
$bc \geq \pi^2 n^2 \bar{c}^2/L^2$.

We now look for $s-$values corresponding to a repeated eigenvalue $\lambda = \frac{-(a-d)}{2\bar{c}}$ (i.e. where $s = \frac{a+d}{2} \pm \sqrt{bc}$, so that $\text{disc}(M) = 0$) with multiplicity 2 and only one eigenvector.  The general solution \eqref{eqn: general} reduces to \begin{align}
\hat{W}(x, s) = &C_1(s)\exp\left[\lambda(s) x\right](I + x(M-\lambda(s) I))v_1 + \notag\\
&C_2(s)\exp\left[\lambda(s) x\right](I + x(M-\lambda(s) I))v_2,
\end{align}
where $v_1, v_2$ are linearly independent solutions of $(M-\lambda I)^2 v = 0$, namely $v_1 = [1, 0]^T$ and $v_2 = [0, 1]^T$. Imposing boundary conditions \eqref{eqn: char_bc} and seeking non-trivial solutions as before yields the equation 

\begin{align}\label{eqn: detA_case2}
    \exp\left[\lambda(s) L\right]&\left[-R_0\left(-\frac{bL}{\bar{c}} - R_L\left[1 + \frac{L}{\bar{c}}\left(\frac{a + d}{2} - s\right)\right]\right) -\right.\notag\\ & \left. \left(1 + \frac{L}{\bar{c}}\left[s - \frac{a + d}{2}\right] - R_L\frac{cL}{\bar{c}}\right) \right]  = 0.
\end{align}
In the case that $R_0 = R_L = R \in \{-1, 1\}$, \eqref{eqn: detA_case2} reduces to 

 \begin{equation}\label{eqn: case2root}
   \exp\left[\lambda(s) L\right]\left[R\frac{L}{\bar{c}}\left(b + c\right) + \frac{L}{\bar{c}}\left(a + d - 2s\right)\right] = 0.
\end{equation}
The second term in \eqref{eqn: case2root} yields the single real root $s_0$.

%If $\text{imag}(\text{disc}(M)) \neq 0$ then additional complex roots will exist that satisfy
%
%bla.
%
%

%If $R = R_0 = -R_L \in \{\{-1\}, \{1\}\}$ then \eqref{eqn: detA} simplifies to 

%\begin{equation}
%    e^{\lambda_1 L}(1 + R v_1^2)(R - v_2^2)   +  %e^{\lambda_2 L}(1 + R v_2^2)(-R + v_1^2) = 0.
%\end{equation}

%One can check that $s_c$ are also roots, along with the real roots
%
%bla
%
%which can be computed numerically. 

For other values of $R_0, R_L$, roots of \eqref{eqn: detA} are more challenging to compute analytically but could be determined numerically.  

\subsection{Case 3: Variable coefficients; $\tilde{B}$ is a diagonal matrix}
We first define the function 
\begin{equation}
    I_{\bar{c}/\delta}(x) = \displaystyle\int_0^x \frac{\delta(y)}{\bar{c}(y)}dy, 
\end{equation}
with $I_{\bar{c}}(x) = \int_0^x \frac{1}{\bar{c}(y)}dy$ corresponding to the time it takes to travel a distance $x$ at speed $\bar{c}(x)$. We can write the general solution \eqref{eqn: general} as 
\begin{align}
    \hat{W}(x, s) = &C_1(s)\exp\left[-I_{\bar{c}/a}(x) + sI_{\bar{c}}(x)\right] v_1 +  \notag \\ &C_2(s) \exp\left[I_{\bar{c}/d}(x) - sI_{\bar{c}}(x)\right] v_2.
\end{align}
Application of boundary conditions and the desire for non-trivial solutions requires again solving $\text{det}(A) = 0$, which yields the equation 
\begin{equation}
    R_0 R_L \exp\left[I_{\bar{c}/d}(L) - sI_{\bar{c}}(L)\right] - \exp\left[-I_{\bar{c}/a}(L) + sI_{\bar{c}}(L)\right] = 0,
\end{equation}
which can be solved for the discrete eigenvalues 
\begin{equation}
    s_n = \frac{\ln(R_0 R_L) + I_{\bar{c}/d}(L) + I_{\bar{c}/a}(L) + 2\pi n i}{2 I_{\bar{c}}(L)}, \quad n \in \mathbb{Z},
\end{equation}
which reduces to \eqref{eqn: s_n} in the case of constant coefficients. 

\subsection{Case 4: Variable coefficients; $\tilde{B}$ is not a diagonal matrix}
For this case it is more challenging to compute an analytic spectrum since the general solution to the differential equation \eqref{eqn: char2} is not easily expressed in closed form \cite{Taylor2011}. However, we include it here as we explore its spectrum numerically in \Cref{sec:discrete}.

\section{SBP-SAT Finite Difference Approximation}
\label{sec:discrete}
Analytic solutions to the governing equations \eqref{eqn: char} - \eqref{eqn: char_ic} can often be difficult to obtain and we turn to numerical methods.  In this work we consider a class of high-order accurate finite difference methods known as summation-by-parts (SBP), with weak enforcement of boundary conditions through the simultaneous-approximation-term (SAT) technique.  As opposed to traditional finite difference methods that ``inject" boundary data by overwriting grid points with the given data, the SAT technique imposes boundary conditions weakly through penalization, so that all grid points approximate both the PDE and the boundary conditions up to a certain level of accuracy.   The combined approach is known as SBP-SAT and provides a provably stable semi-discretization in space \cite{DELREYFERNANDEZ2014171, SVARD201417, CARPENTER1994220}. The spatial derivative operators involve centered difference approximations within the domain interior, with a special transition to one-sided differences near domain boundaries. We use the diagonal norm SBP-SAT operators from \cite{strand_summation_1994} that have a formal order of accuracy given by $p = 2; 3; 4$ and $5$, with interior order of accuracy $2p - 2$ and boundary accuracy of $p - 1$.

\subsection{Discrete Operators}

The spatial domain $0 \leq x \leq L$ is discreteized with $N+1$ evenly spaced grid points $x_i = 0 + ih, i=0, \dots, N$ with grid spacing $h = 1/N$. A function $u$ projected onto the computational grid is denoted by $\mm{u} = [u_0, u_1, \dots, u_N]^T$ and is often taken to be the interpolant of $u$ at the grid points. The grid basis vector $\vv{e}_j$ is defined to be a vector with value 1 at grid point $j$ and 0 otherwise, which allows us to extract the jth component, namely $u_j = \vv{e}_j^T\vv{u}$.  
 
 We say that matrix $\mm{D}$ is an SBP approximation to the first derivative operator $\partial /\partial
  x$ if it can be decomposed as $\mm{H}\mm{D} = \mm{Q}$ with $\mm{H}$ being symmetric positive-definite and $\mm{Q}$ satisfying $\vv{u}^{T}(\mm{Q} +
  \mm{Q}^{T})\vv{v} = u_{N}v_{N} - u_{0}v_{0}$. Here, $\mm{D}$ is the standard central finite difference operator in the interior which transitions to one-sided at boundaries and $\mm{H}$ is a diagonal, positive-definite matrix allowing us to define a discrete inner-product and norm, namely
\begin{equation}
    (\mm{U}, \mm{W})_\mm{H} = \mm{U}^T \mm{H}_2 \mm{W}, \quad ||\mm{U}||^2_\mm{H} = (\mm{U}, \mm{U})_\mm{H},
    \end{equation}
where $\mm{H}_2 = (\mm{I}_2 \otimes \mm{H})$ and $\mm{U} = [\mm{u}_1^T \quad \mm{u}_2^T]^T, \mm{W} = [\mm{w}_1^T \quad \mm{w}_2^T]^T$.
\subsection{Discretized Governing Equations}
Equation \eqref{eqn: char} and boundary conditions \eqref{eqn: char_bc} are discretized with the skew-symmetric SBP-SAT method \cite{NORDSTROM2017451, FISHER2013353} as 
   \begin{align}\label{eqn: discrete1}
    \mm{W}_t + \mm{\Lambda}_D\mm{W} = & \tilde{\mm{B}}\mm{W}  + \alpha_0\mm{H}_2^{-1}(\mm{w}_{2, 0} - R_0 \mm{w}_{1, 0} - h_0(t))(E_2 \otimes e_0)  + \notag \\& \quad \alpha_L\mm{H}_2^{-1}(\mm{w}_{1,N} - R_L\mm{w}_{2, N} - h_L(t))(E_1 \otimes e_N) + \tilde{\mm{F}},
\end{align}
where \begin{equation}
    \mm{\Lambda}_D = 0.5(\mm{\Lambda}  \mm{D}_2 + \mm{D}_2 \mm{\Lambda}) - 0.5 \mm{M}_{\bar{c}},
\end{equation}
for $\mm{\Lambda} = \text{diag}(-\mm{\bar{c}}, \mm{\bar{c}})$ and $\mm{D}_2 = (\mm{I}_2 \otimes \mm{D})$.  Here 
\begin{equation}\mm{M}_{\bar{c}}=\begin{pmatrix}-\text{diag}(\mm{D}\mm{\bar{c}})&\mm{0}\\\mm{0}&\text{diag}(\mm{D}\mm{\bar{c}})\end{pmatrix}, \quad \tilde{\mm{B}} = \begin{pmatrix}\text{diag}(\mm{a})&\text{diag}(\mm{b})\\\text{diag}(\mm{c})&\text{diag}(\mm{d})\end{pmatrix},\end{equation}
which reduces to 
$\mm{\Lambda}_D = (\Lambda \otimes \mm{D}), \tilde{\mm{B}} = (\tilde{B} \otimes \mm{I})$ for constant coefficients. 
$\mm{I}$ and $\mm{I}_2$ are the $N \times N$ and $2 \times 2$ identity operators, respectively. Here $E_j$ is a $2 \times 1$ vector of zeros with a 1 in the $j^{th}$ position and $\alpha_0, \alpha_L$ are scalar penalty parameters. A discrete energy estimate determines values of the penalty parameters that yield a stable scheme, which may be obtained by taking $h_0 = h_L = \tilde{F}= 0$ (not required but eases the computation). To aid the following analysis we rewrite \eqref{eqn: discrete1} as 
   \begin{align}\label{eqn: discrete2}
    \mm{W}_t + \mm{\Lambda}_D\mm{W} = &\tilde{\mm{B}}\mm{W} + \mm{P}\mm{W},
\end{align}
where
\begin{align}
       \mm{P} &= \alpha_0\mm{H}_2^{-1}[(E_2 E_2^T \otimes e_0 e_0^T) - R_0 (E_2 E_1^T \otimes e_0 e_0^T)] + \notag \\ &\alpha_L\mm{H}_2^{-1}[(E_1 E_1^T \otimes e_N e_N^T) - R_L (E_1 E_2^T \otimes e_N e_N^T)].
       \end{align}
 Multiplying \eqref{eqn: discrete2} by $\mm{W}^T \mm{H}_2$ and adding the transpose yields the estimate
 \begin{align}\label{eqn: discrete_energy1}
    \dot{\mm{E}_h} &= -\bar{c}_0w_{1, 0}^2 + \bar{c}_N w_{1, N}^2 + \bar{c}_0 w_{2, 0}^2 -\bar{c}_N w_{2, N}^2 + 0.5(\mm{W}, \mm{M}_{\bar{c}}\mm{W}) + 0.5(\mm{M}_{\bar{c}}\mm{W},\mm{W}) + \notag\\ &\quad  (\mm{W}, \tilde{\mm{B}}\mm{W}) + (\tilde{\mm{B}}\mm{W}, \mm{W})+  (\mm{W}, \mm{P}\mm{W}) + (\mm{P}\mm{W},\mm{W}) \notag\\
    & = -\bar{c}_0w_{1, 0}^2 + \bar{c}_N w_{1, N}^2 +\bar{c}_0 w_{2, 0}^2 -\bar{c}_N w_{2, N}^2 + 0.5(\mm{W}, \mm{M}_{\bar{c}}\mm{W}) + 0.5(\mm{M}_{\bar{c}}\mm{W},\mm{W}) + \notag \\
    &\quad  (\mm{W}, \tilde{\mm{B}}\mm{W}) + (\tilde{\mm{B}}\mm{W}, \mm{W}) + \notag \\
    &\quad 2\alpha_0 w_{2,0}(w_{2,0} - R_0 w_{1,0}) + 2\alpha_L w_{1,N}(w_{1, N} - R_L w_{2, N}),
\end{align}
where $\mm{E}_h = \mm{W}^T \mm{H}_2 \mm{W}$.  Taking $\alpha_0 = -\bar{c}_0, \quad \alpha_L = -\bar{c}_N$, estimate \eqref{eqn: discrete_energy1} reduces to 
  \begin{align}\label{eqn: discrete_energy2}
    \dot{\mm{E}_h} = -&\bar{c}_0[w^2_{1, 0}(1-R_0^2) + (w_{2,0} - R_0 w_{1, 0})^2] + \notag \\ -&\bar{c}_N[w_{2,N}^2(1-R_N^2) + (w_{1, N} - R_Lw_{2, N})^2]  + \notag \\&0.5(\mm{W}, \mm{M}_{\bar{c}}\mm{W}) + 0.5(\mm{M}_{\bar{c}}\mm{W},\mm{W}) + (\mm{W}, \tilde{\mm{B}}\mm{W}) + (\tilde{\mm{B}}\mm{W}, \mm{W}) \notag \\
   & \leq (||\mm{M}_{\bar{c}}||_\mm{H} + 2||\tilde{\mm{B}}||_\mm{H})\mm{E}_h,
\end{align}
 and thus a stable discretization, where the rate of change of the discrete energy mimics that of continuous case as in \eqref{eqn: energydot}. 

We note that the discrete energy estimate \eqref{eqn: discrete_energy2}, places bounds on the solution at any final time $T$, for the case of zero boundary and source data. This means that the numerical solution converges as the grid is refined but only on a compact domain in both space and time. Because we are interested in how a discretization approximates analytic growth/decay rates it is important to underscore that such stability estimates do not guarantee the temporal evolution of the discrete system matches that of the of the continuous (which is ensured by a stronger notion known as strict stability, see \cite{SVARD201417}), which we will discuss in more detail in \Cref{sec:discrete_spec}.

\section{Analytic Solutions and Numerical Verification}
\label{sec:verification}
To verify convergence of our numerical methods, analytic solutions are desirable. 
When $\tilde{B}$ is a diagonal matrix and model data (e.g. $\tilde{F}, h_0$ etc.) are mostly zero, analytic solutions may be obtained via Laplace transform methods as in \cite{Erickson2019} without computations becoming too untractable. For more general cases of $\tilde{B}$ and non-zero model data however, we may obtain analytic solutions using the method of manufactured solutions \cite{Roache}.  

To illustrate the construction of an analytic solution via Laplace methods, we consider case 3 with $b(x) = c(x) = 0$, $h_0 = h_L = F = 0$, and initial data $f_1(x) = 0, f_2(x) \neq 0$. The solution to \eqref{eqn: char}-\eqref{eqn: char_ic} is thus
\begin{equation}
    \hat{W} = \exp\left(\int_0^x M(\xi, s) d\xi\right)C(s) + G(x, s),
\end{equation}
where \begin{equation}
   G(x, s) = [0, \quad G_2(x, s)]^T, 
\end{equation}
for 
\begin{equation}
  G_2(x, s) = \displaystyle\int_0^x \frac{f_2(\xi)}{\bar{c}(\xi)} \exp\left[I_{\bar{c}/d}(x) - I_{\bar{c}/d}(\xi) - s(I_{\bar{c}}(x) - I_{\bar{c}}(\xi))\right] d\xi,
\end{equation}
% \bar{c}/d
or in component form
\begin{subequations}\label{eqn: comp}
    \begin{align}
        \hat{w}_1(x, s) &= \exp\left[-I_{\bar{c}/a}(x) + sI_{\bar{c}}(x)\right]C_1(s),\\
        \hat{w}_1(x, s) &= \exp\left[I_{\bar{c}/d}(x) - sI_{\bar{c}}(x)\right]C_2(s) + G_2(x, s).
    \end{align}
\end{subequations}
Using similar techniques as in \cite{Erickson2019}, we observe that 
\begin{equation}
    G_2(x, s) = \exp\left[I_{\bar{c}/d(x)}\right]\mathcal{L}[f_2(\xi)H(\xi)\exp\left[-I_{\bar{c}/d}(\xi)\right]],
\end{equation}
where $H$ is the Heaviside function and $\xi$ is the characteristic variable defined implicitly by 
\begin{equation}
    t = I_{\bar{c}}(x) - I_{\bar{c}}(\xi).
\end{equation}
We note that closed forms solutions therefore require that $I_{\bar{c}}(x)$ be invertible, which allows us to solve for $\xi$, namely 
\begin{equation}\label{eqn: xi}
    \xi(x, t) = I_{\bar{c}}^{-1}(I_{\bar{c}}(x) - t),
\end{equation}
as in \cite{Erickson2019}. Applying boundary conditions allows us to solve for the coefficients, namely, we have
\begin{align}
    C_1(s) &= R_L R(s) G_2(L, s),\\
    C_2(s) &= R_0 R_L R(s) G_2(L, s),
\end{align}
where 

\begin{align}\label{eqn: R}
    R(s) = &\frac{1}{\exp\left[-I_{\bar{c}/a}(L) + sI_{\bar{c}}(L)\right] - R_0 R_L \exp\left[I_{\bar{c}/d}(L) - s I_{\bar{c}}(L)\right]} \\
= &\exp\left[I_{\bar{c}/a}(L) - s I_{\bar{c}}(L)\right] \cdot \notag \\&\displaystyle\sum_{n = 0}^\infty (R_0 R_L)^n \exp\left[n(I_{\bar{c}/d}(L) + I_{\bar{c}/a}(L) - 2sI_{\bar{c}}(L))\right],
\end{align}
where the geometric series in \eqref{eqn: R} converges for \begin{equation}
    \text{Re}(s) > \frac{\ln(R_0 R_L) + I_{\bar{c}/a}(L) + I_{\bar{c}/d}(L)}{2I_{\bar{c}}(L)}.
\end{equation}
Thus we have solutions in Laplace space given by 
\begin{subequations}
\begin{align}
    \hat{w}_1(x, s) = &\displaystyle\sum_{n = 0}^\infty R_L (R_0 R_L)^n \notag \\
    &\exp\left[y_n(x)-s \alpha_n(x)\right] \mathcal{L}[f_2(\xi)H(\xi)\exp\left[-I_{\bar{c}/d}(\xi)\right]],\\
     \hat{w}_2(x, s) = &G_2(x, s) + \displaystyle\sum_{n = 0}^\infty (R_0 R_L)^{n+1} \cdot \notag \\ &\exp\left[z_n(x)-s \beta_n(x)\right] \mathcal{L}[f_2(\xi)H(\xi)\exp\left[-I_{\bar{c}/d}(\xi)\right]] ,
     \end{align}
\end{subequations}
where 
\begin{align}
y_n(x) &= (-I_{\bar{c}/a}(x) + I_{\bar{c}/a}(L) + n[I_{\bar{c}/d}(L) + I_{\bar{c}/a}(L)] + I_{\bar{c}/d}(L),\\
z_n(x) & = (I_{\bar{c}/d}(x) + I_{\bar{c}/a}(L) + n[I_{\bar{c}/d}(L) + I_{\bar{c}/a}(L)] + I_{\bar{c}/d}(L),\\
\alpha_n(x) &= -I_{\bar{c}}(x) + (2n+1)I_{\bar{c}}(L),\\
\beta_n(x) & = I_{\bar{c}}(x) + (2n+1)I_{\bar{c}}(L)
\end{align}
which reduces to 
\begin{align}
y_n(x) &= (-ax + aL + nL(a+d) + dL)/\bar{c},\\
z_n(x) & = (dx + aL + nL(a+d) + dL)/\bar{c},\\
\alpha_n(x) &= ((2n+1)L - x)/\bar{c},\\
\beta_n(x) & = ((2n+1)L + x)/\bar{c}
\end{align}
for constant coefficients.
Applying the inverse Laplace transform to \eqref{eqn: comp} then yields the following solutions
\begin{subequations}\label{eqn: final4}
    \begin{align}
       w_1(x, t) = &\displaystyle\sum_{n = 0}^\infty R_L(R_0 R_L)^n \cdot \notag \\
       &\exp\left[y_n(x)-I_{\bar{c}/d}(\xi_n)\right] f_2(\xi_n) H(\xi_n) H(t-\alpha_n(x)),\\
       w_2(x, t) = & \exp\left[I_{\bar{c}/d}(x)\right]f_2(\xi)H(\xi)\exp\left[I_{\bar{c}/d}(\xi)\right] + \notag\\
       &\displaystyle\sum_{n = 0}^\infty (R_0 R_L)^{n+1} \cdot \notag \\ &\exp\left[z_n(x)-I_{\bar{c}/d}(\gamma_n)\right] f_2(\gamma_n) H(\gamma_n) H(t-\beta_n(x)),
    \end{align}
\end{subequations}
where 
\begin{subequations}
\begin{align}
    \xi_n(x, t) &= \xi(L, t - \alpha_n(x)),\\
    \gamma_n(x, t) &= \xi(L, t - \beta_n(x)),
\end{align}
\end{subequations}
where $\xi(x, t)$ is defined by \eqref{eqn: xi} and the argument $(x, t)$ is dropped in \eqref{eqn: final4} for notational ease. 
The solutions in \eqref{eqn: final4} illustrate how $\bar{c}(x)$ determines the speed at which the initial data $f_2(x)$ is translated back and forth across the domain. Parameter $a(x)$ contributes to the spatial structure of the solution, whereas $d(x)$ (but not $a(x)$) contributes to exponential growth or decay in time, depending on its sign. However, if $f_1(x) \neq 0$ then $a$ would also contribute to temporal exponential growth/decay. 

For more complex solutions (for example where $h_0, h_L, F \neq 0$),  one may manufacture an analytic solution. For example here we assume the analytic solution is (denoted with asterisks)
\begin{subequations}\label{eqn: manuf}
    \begin{align}
        w^*_1(x, t) &= \sin(\pi x - t)\\
        w^*_2(x, t) & = \cos(2\pi x + t).
    \end{align}
\end{subequations}
Thus for any values of $\bar{c}$ and matrix $\tilde{B}$, the corresponding source term is $\tilde{F}(x, t) = [\tilde{F}_1(x, t), \quad \tilde{F}_2(x, t)]^T$ where
\begin{align}
    \tilde{F}_1(x, t) = &-(1 + \bar{c}(x) \pi)\cos(\pi x - t) + \notag \\
    &-a(x) \sin(\pi x - t) - b(x)\cos(2\pi x + t),\\ 
    \tilde{F}_2(x, t) = &-(1 + 2\bar{c}(x)\pi)\sin(2\pi x +t) + \notag \\
    &-c(x) \sin(\pi x - t) - d(x)\cos(2\pi x + t),
\end{align}
and the boundary data is
\begin{subequations}
\begin{align}
       h_0(t) &= \cos(t) + R_0 \sin(t),\\
       h_L(t) &= \sin(\pi L - t) - R_L \cos(2\pi L + t),
\end{align}
\end{subequations}
for any length $L$ and reflection coefficients. 

\begin{figure*}[t!]
\centering
        \includegraphics[height=3in]{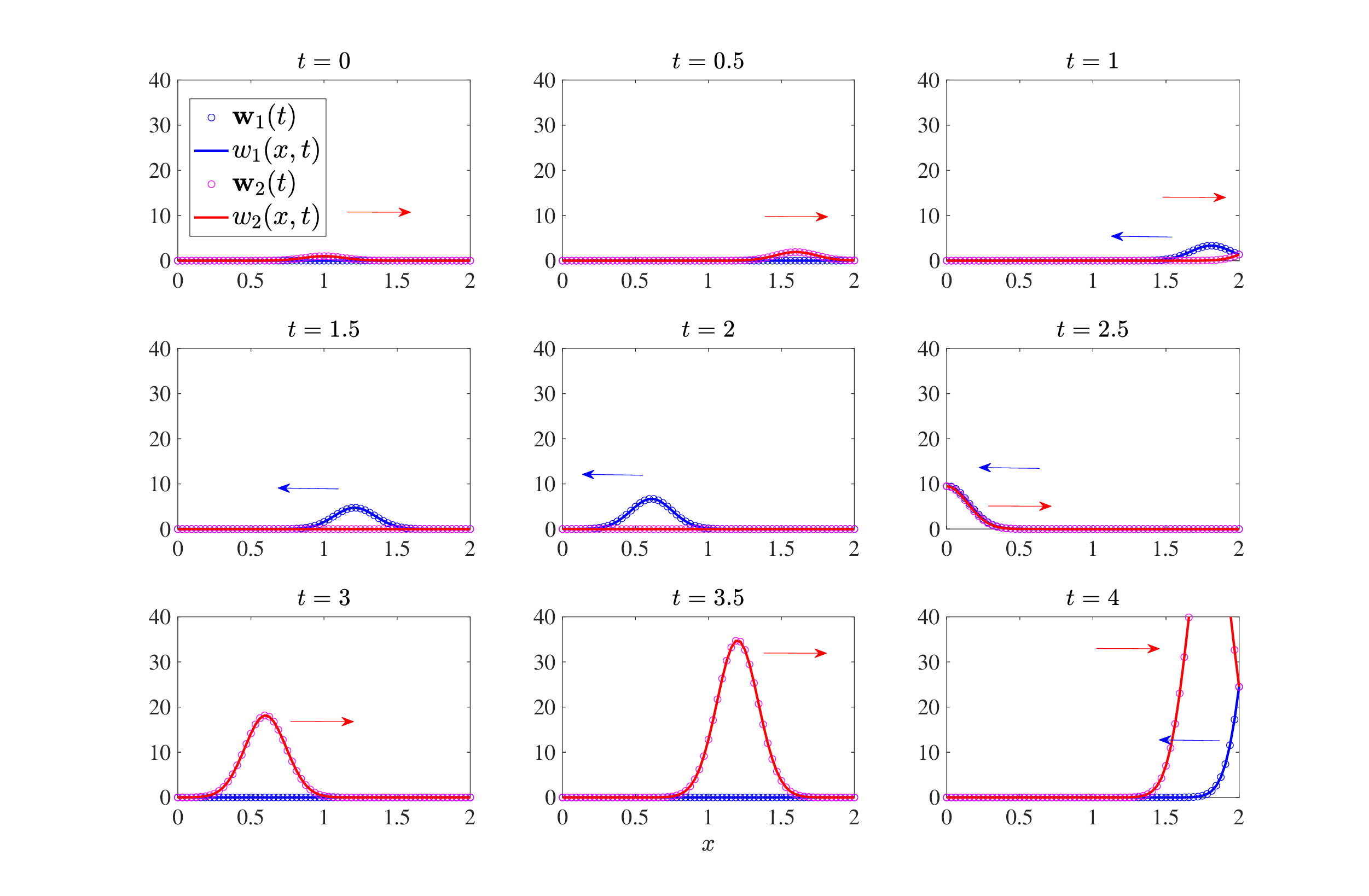}
    \caption{Temporal evolution of the analytic and numerical approximation to the solution corresponding to Case 1b using $p = 2$, $N = 2^{10}$, where exponential growth is present in the physical problem. Initial condition is a Gaussian pulse on $w_2$ (red), centered at zero, which initially propagates to the right (direction of pulse propagation denoted with arrows and colors correspond to the solution component), is reflected to $w_1$ (blue) which then propagates to the left etc. Time denoted above subfigures.   }
    \label{fig: analy_1b}
\end{figure*}

For the numerical verification we take $R_0 = R_L = 1$, $L = 2$ and consider two variants of all 4 cases, listed in \Cref{table1}. Note that we limit our investigation of variable coefficients to smooth, linear functions; the effects of nonlinear and/or non-smooth coefficients can be found in \cite{Nord2006, Erickson2019}.
\begin{table}[t]
\centering
\begin{tabular}{l l l l l l}%% Table column specifiers
%% Tabular cells are separated by &
  Case & $a(x)$ & $b(x)$ & $c(x)$ & d(x) & $\bar{c}(x)$ \\ %% A tabular row ends with \\
  \hline
  1a & -0.7 & 0 & 0 & -1.3 & 1.2 \\
1b & 0.7 & 0 & 0 & 1.3& 1.2 \\
2a & 1 & -3 & 3 & -5 & 1.2\\
2b & 2 & 3 & 2 & 4 & 1.2\\
3a & $-0.7 - x$ & 0 & 0 & $-1.3 + 0.1x$ & $1.2 + 0.5x$\\
3b &$ 0.7 + x$ & 0 & 0 & $1.3 + 0.1x$& $1.2 + 0.5x$ \\
4a & $x$ & $-3x$ & $3x$ & $-5x$ & $1.2 + 0.5x$\\
4b & $2 + x$ & $3x$ & $2x$ & $-1 + 4x$ & $1.2 + 0.5x$\\
\end{tabular}
% Use \caption command for table caption and label.
\caption{Values for the variable coefficients for Cases 1-4. }\label{table1}
\end{table}
For Cases 1 and 3 we verify with the analytic solution \eqref{eqn: final4} using $f_2(x) = \exp\left[-(x - L/2)^2/\sigma^2\right]$ with $\sigma = 0.2$. The temporal evolution of the analytic and numerical solutions for case 1b (where physical growth is present) is illustrated in \Cref{fig: analy_1b}. For Cases 2 and 4 we use the manufactured solution \eqref{eqn: manuf}.  Initial conditions, boundary and source terms are set by the exact solutions and we simulate over the temporal interval $0 \leq T$. For case 1 we set $T = 3$ and for cases 2-4 we take $T = 0.1$. 

\begin{figure*}[t!]
\centering
    \begin{subfigure}[t]{0.5\textwidth}
        \centering
        \includegraphics[height=1.7in]{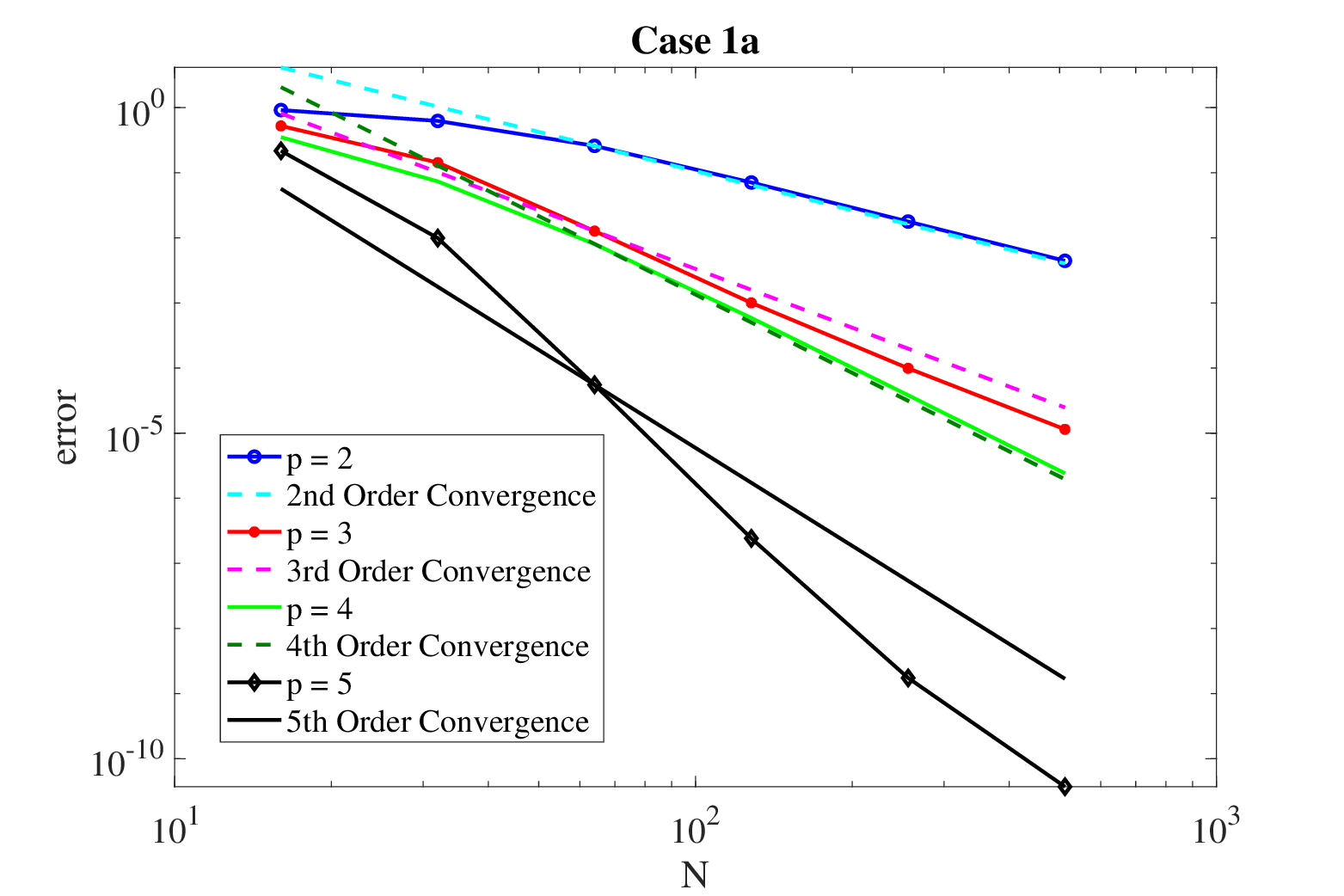}
        \caption{}
    \end{subfigure}%
    ~ \hspace{-5mm}
     \begin{subfigure}[t]{0.5\textwidth}
        \centering
        \includegraphics[height=1.7in]{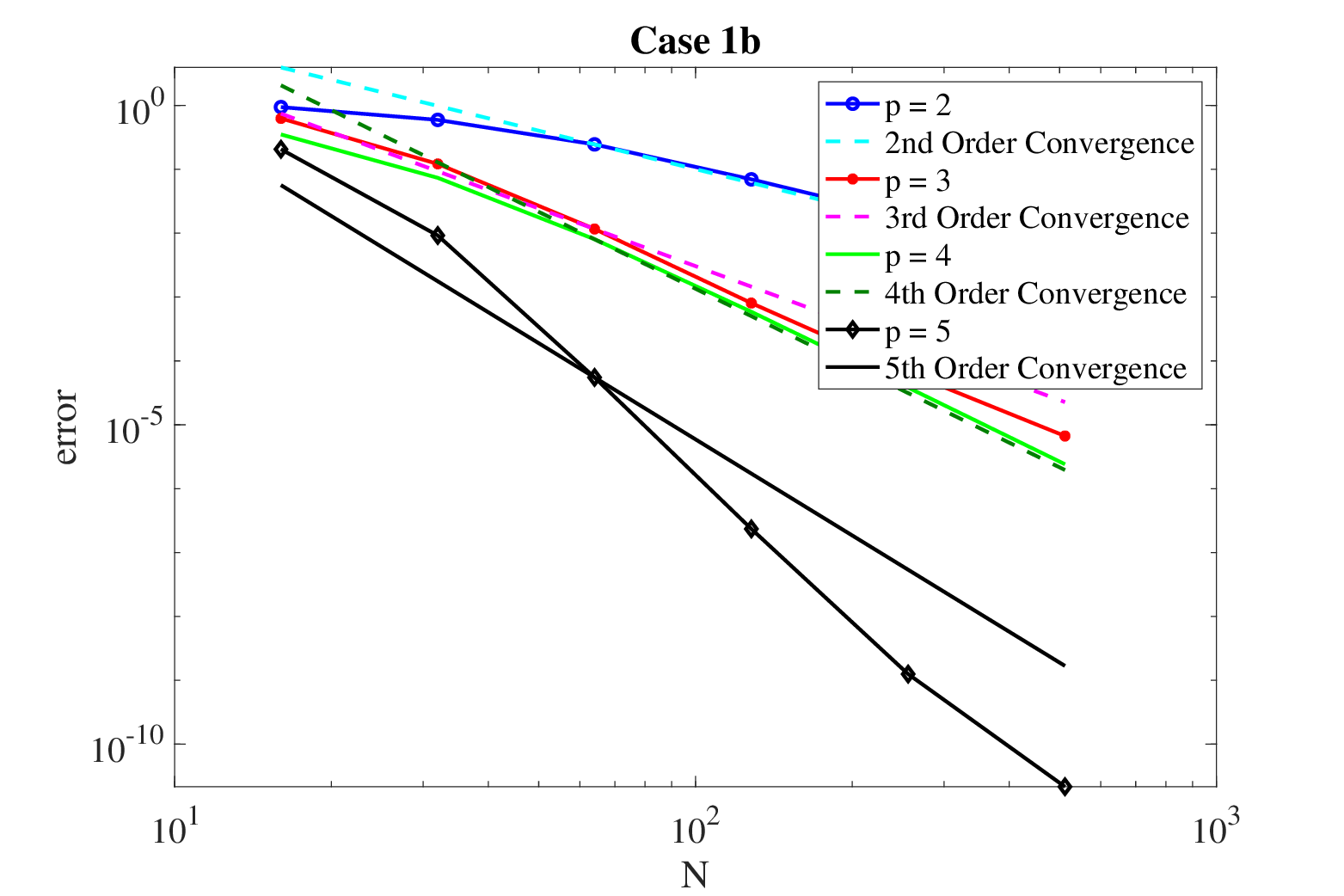}
        \caption{}
    \end{subfigure}
    \caption{Convergence rates for SBP operators with order of accuracy $p = 2; 3; 4; 5$ for (a) case 1a and (b) case 1b. Relative
error computed at $T = 3$, with convergence achieved at the expected rates.}
    \label{fig: conv1}
\end{figure*}

\begin{figure*}[t!]
\centering
    \begin{subfigure}[t]{0.5\textwidth}
        \centering
        \includegraphics[height=1.7in]{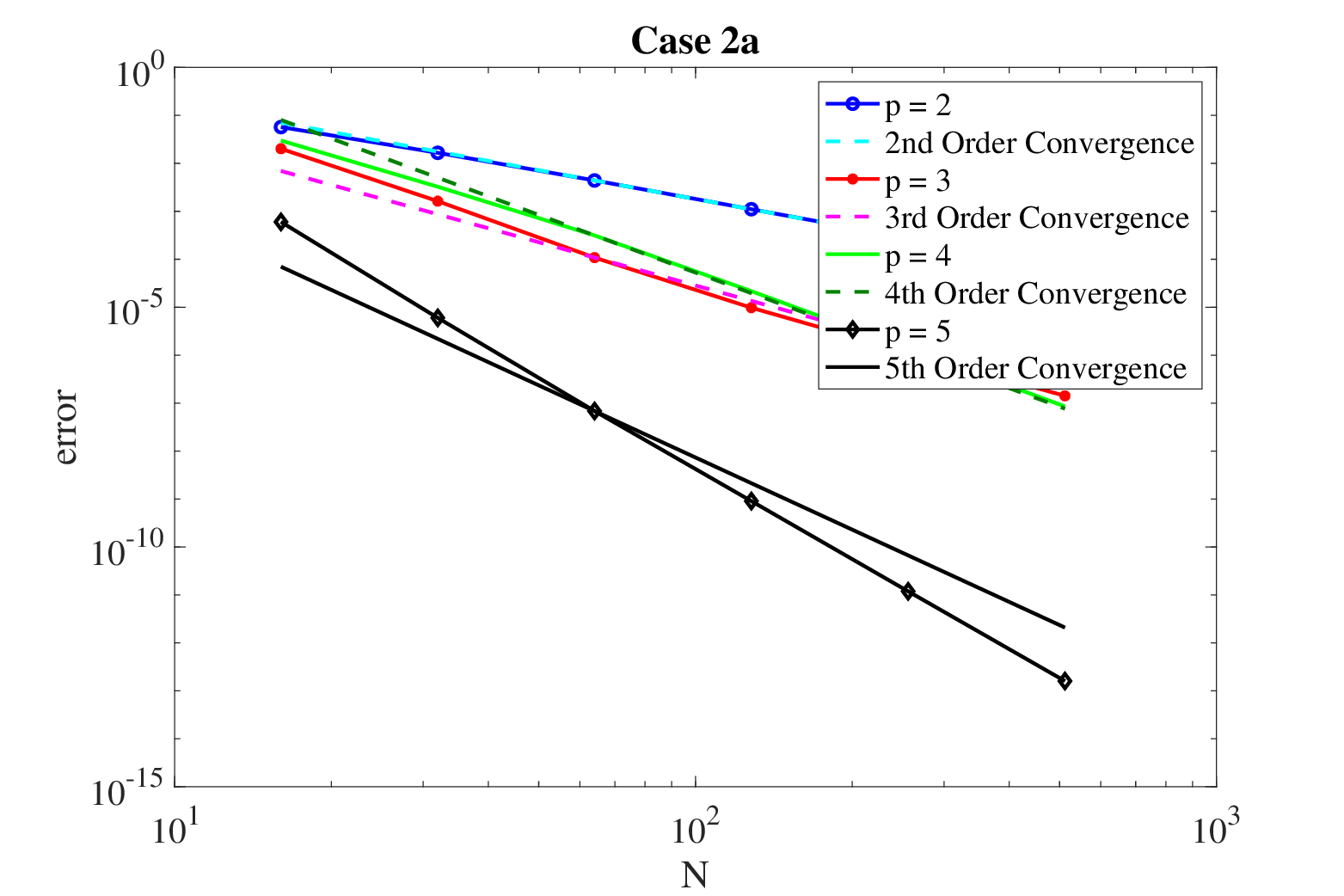}
        \caption{}
    \end{subfigure}%
    ~ \hspace{-5mm}
     \begin{subfigure}[t]{0.5\textwidth}
        \centering
        \includegraphics[height=1.7in]{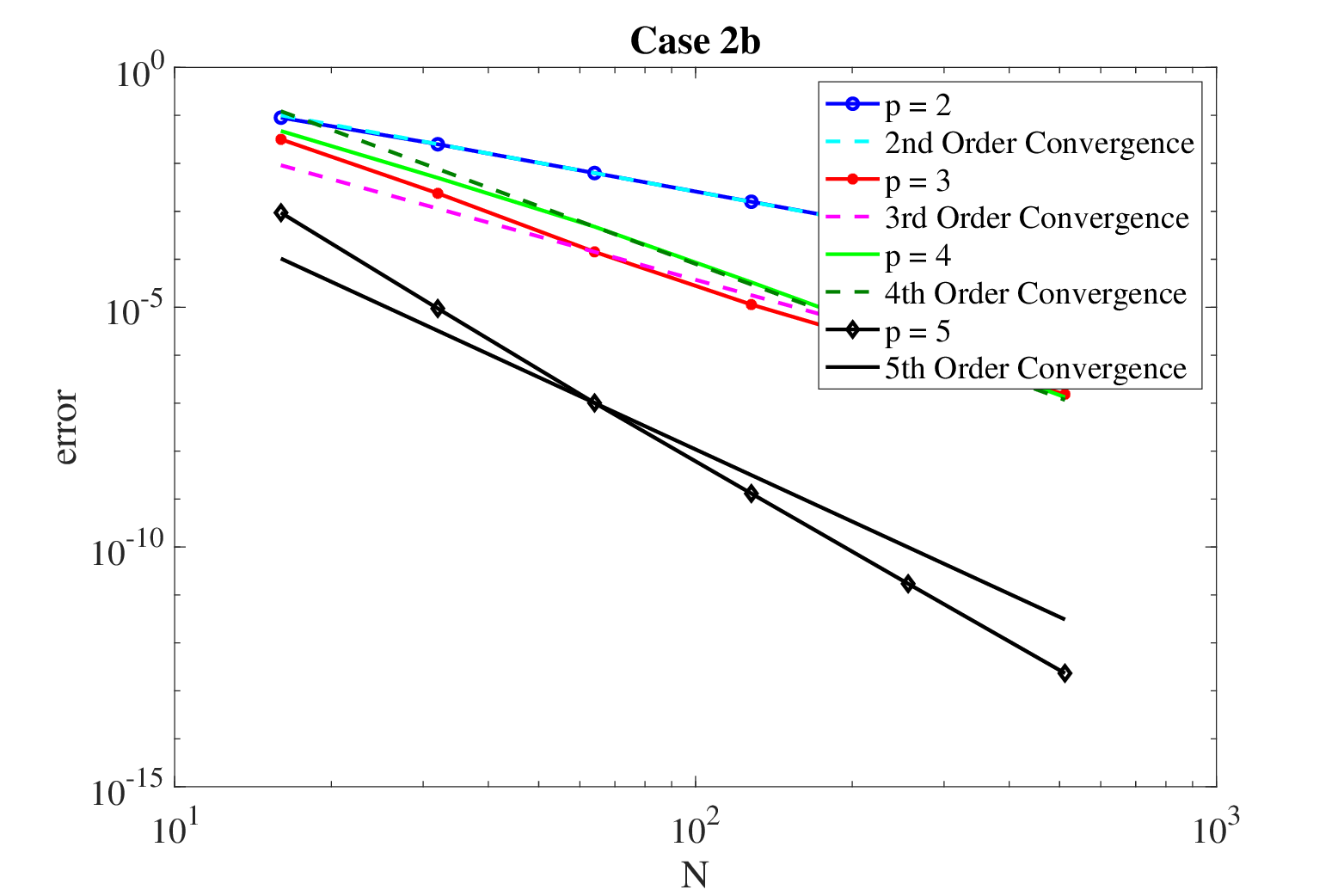}
        \caption{}
    \end{subfigure}
    \caption{Convergence rates for SBP operators with order of accuracy $p = 2; 3; 4; 5$ for (a) case 2a and (b) case 2b. Relative
error computed at $T = 0.1$, with convergence achieved at the expected rates.}
    \label{fig: conv2}
\end{figure*}

\begin{figure*}[t!]
\centering
    \begin{subfigure}[t]{0.5\textwidth}
        \centering
        \includegraphics[height=1.7in]{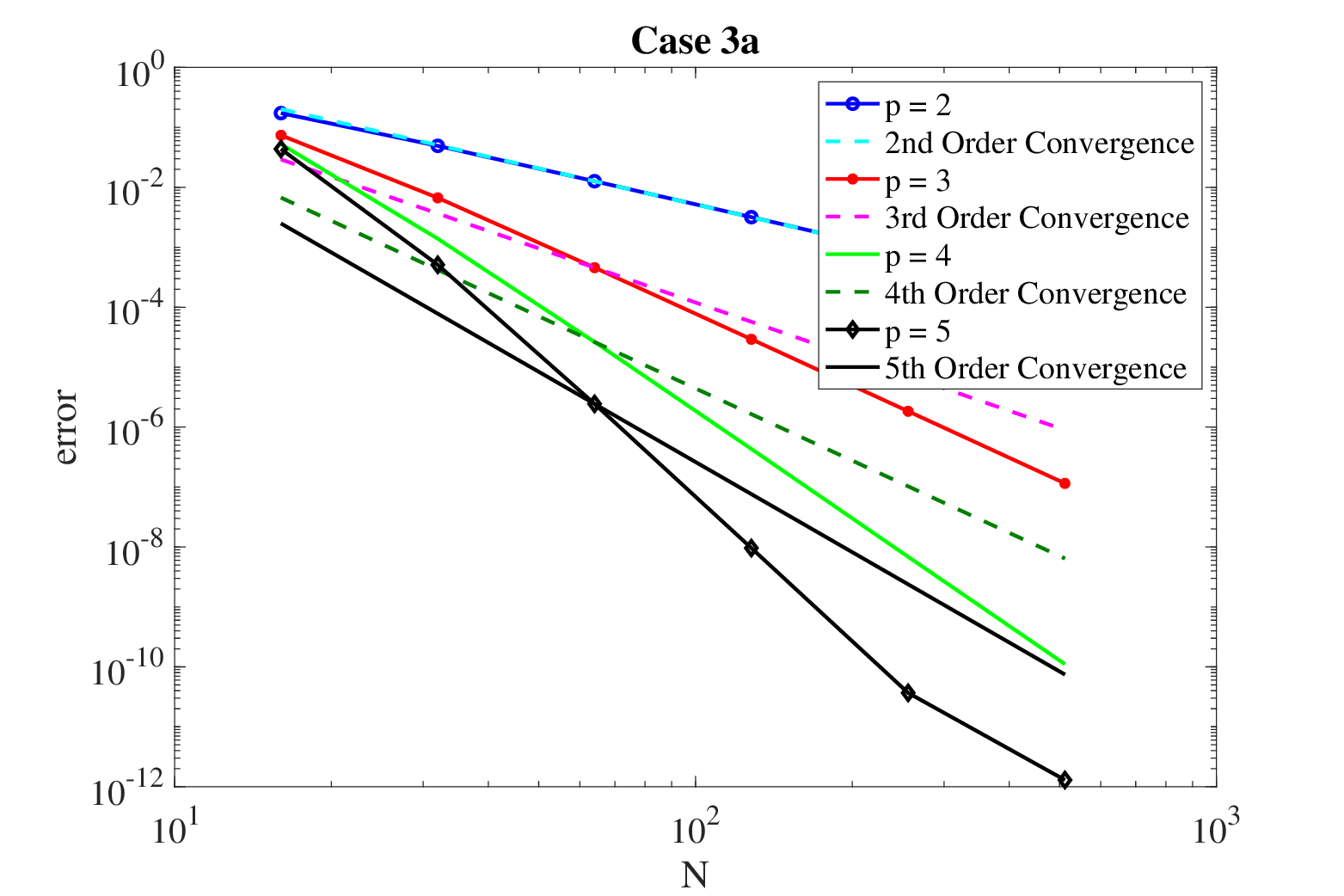}
        \caption{}
    \end{subfigure}%
    ~ \hspace{-5mm}
     \begin{subfigure}[t]{0.5\textwidth}
        \centering
        \includegraphics[height=1.7in]{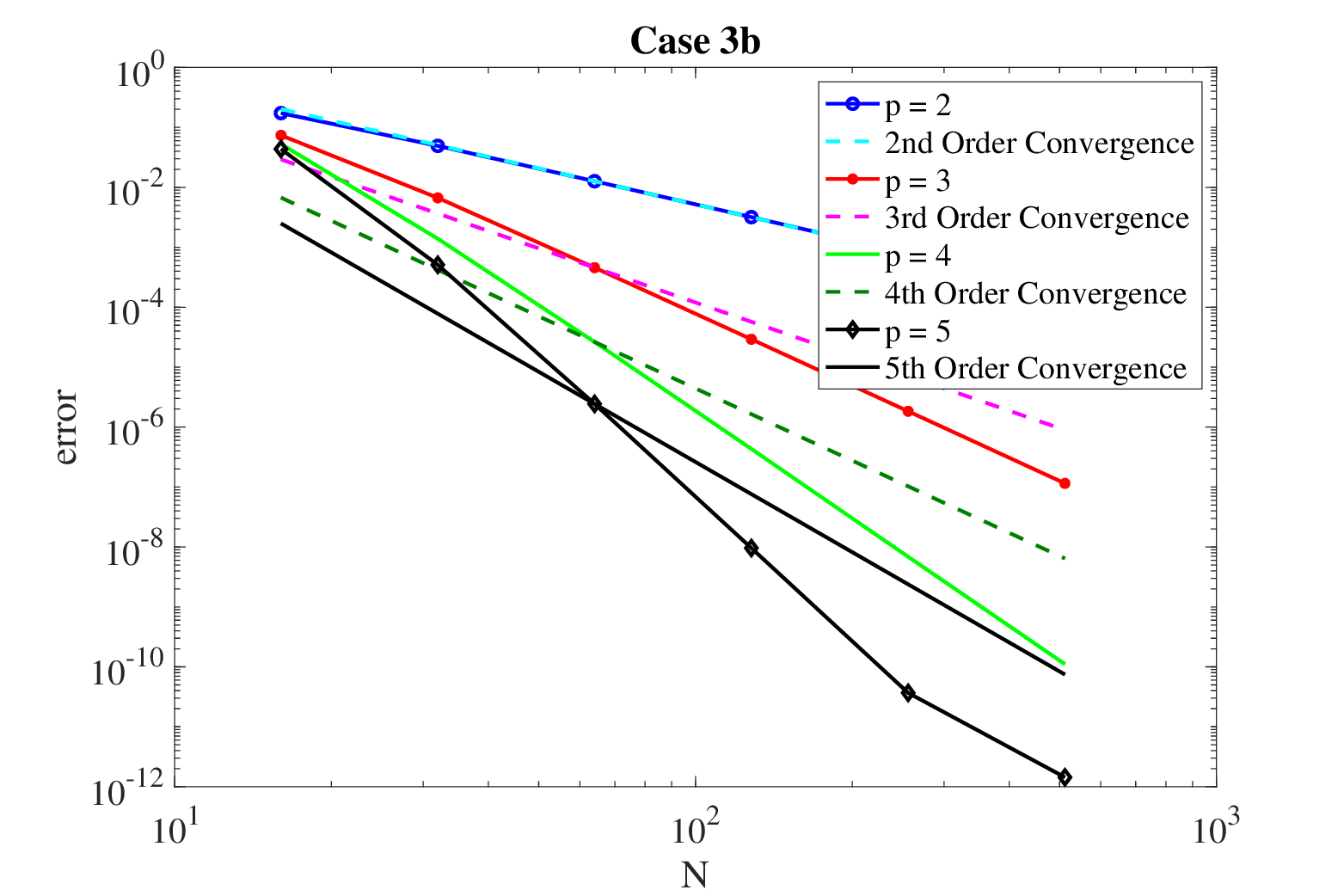}
        \caption{}
    \end{subfigure}
    \caption{Convergence rates for SBP operators with order of accuracy $p = 2; 3; 4; 5$ for (a) case 3a and (b) case 3b. Relative
error computed at $T = 0.1$, with convergence achieved at the expected rates.}
    \label{fig: conv3}
\end{figure*}

\begin{figure*}[t!]
\centering
    \begin{subfigure}[t]{0.5\textwidth}
        \centering
        \includegraphics[height=1.7in]{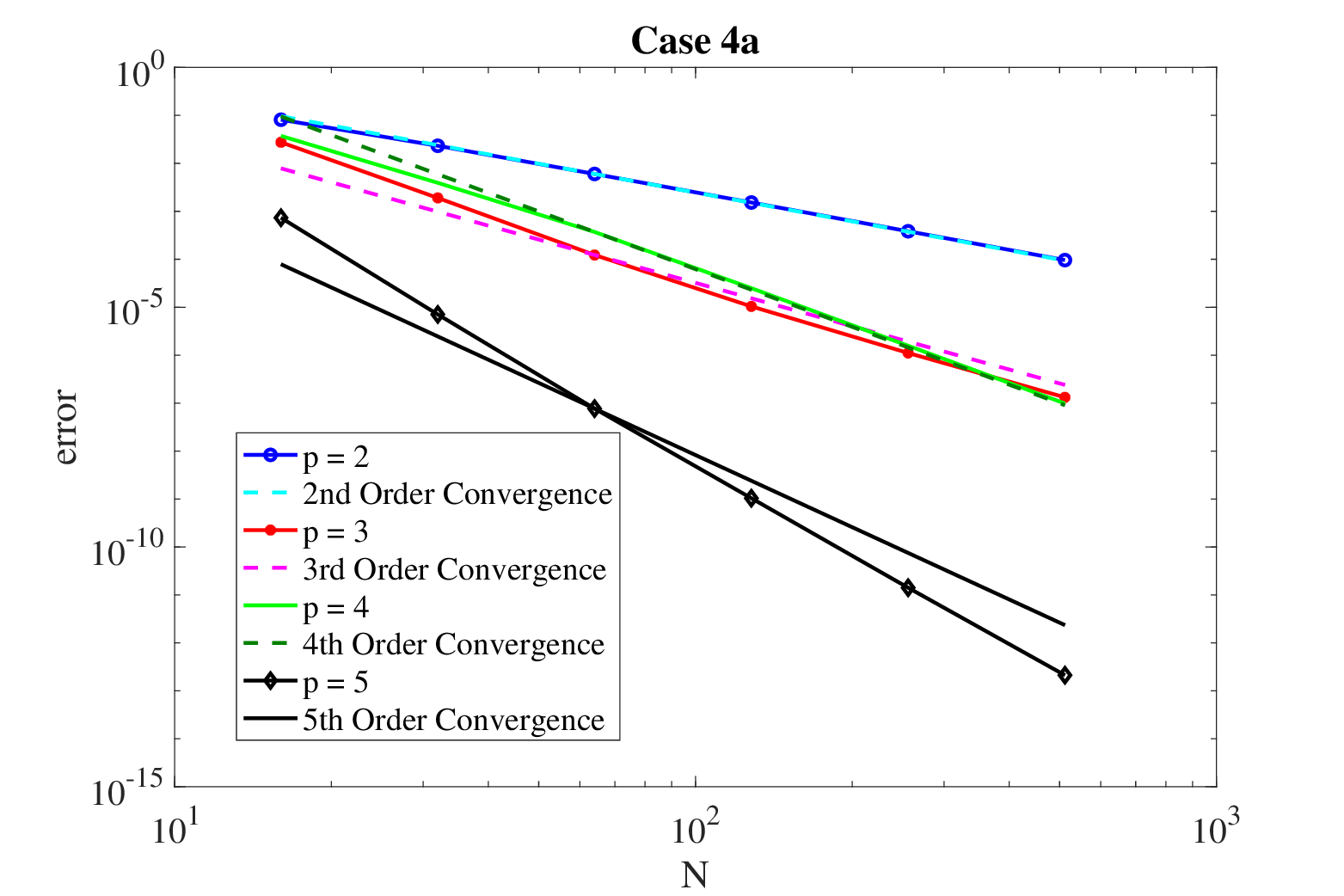}
        \caption{}
    \end{subfigure}%
    ~ \hspace{-5mm}
     \begin{subfigure}[t]{0.5\textwidth}
        \centering
        \includegraphics[height=1.7in]{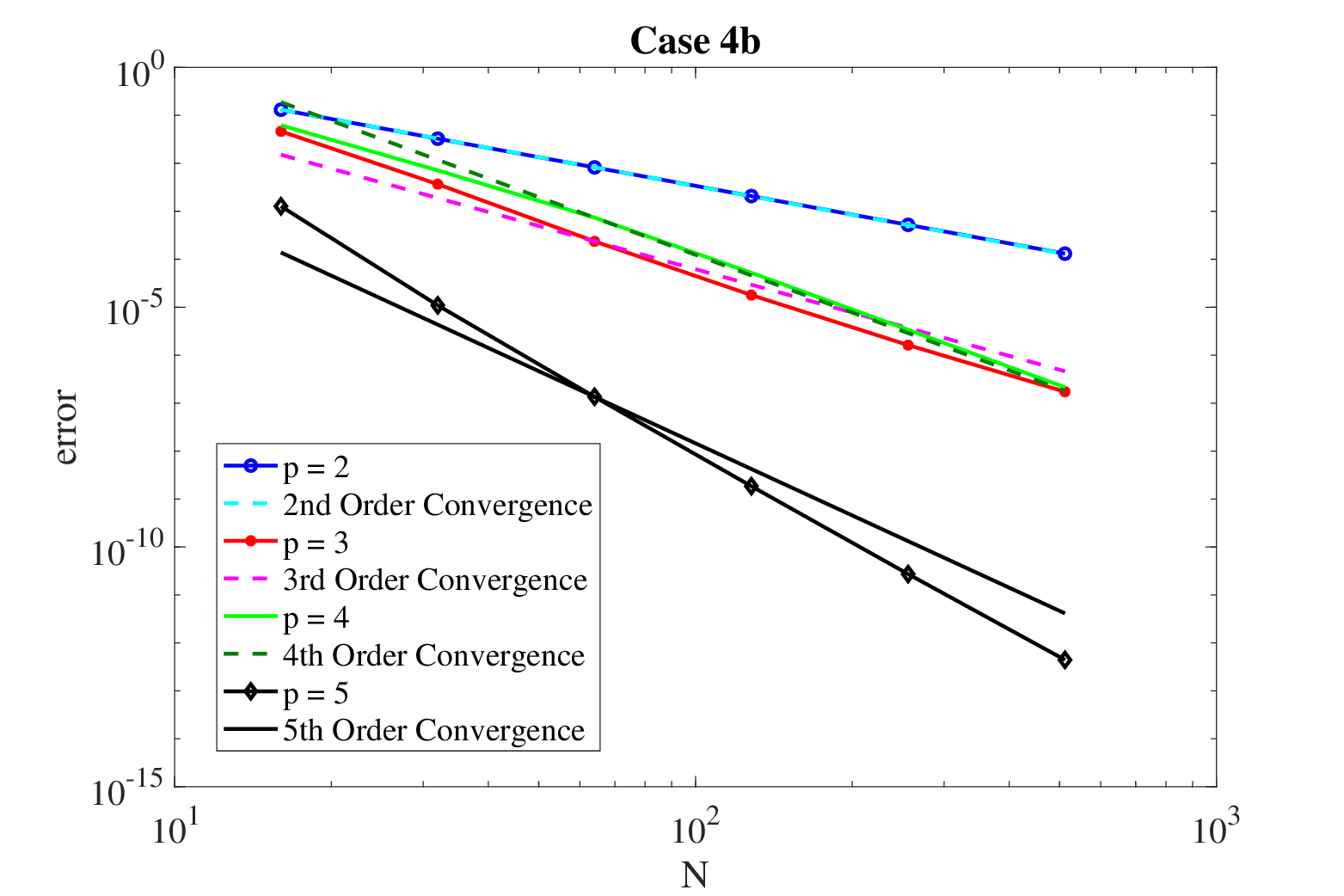}
        \caption{}
    \end{subfigure}
    \caption{Convergence rates for SBP operators with order of accuracy $p = 2; 3; 4; 5$ for (a) case 4a and (b) case 4b. Relative
error computed at $T = 0.1$, with convergence achieved at the expected rates.}
    \label{fig: conv4}
\end{figure*}

We denote the total error in the discrete $\mm{H}-$norm by 
\begin{equation}
    E(t) = ||W^*(\cdot, t) - \mm{W}(\cdot, t)||_\mm{H},
\end{equation}
where $W^* = [w^{*, T}_1 \quad w^{*, T}_2]^T$ is the stacked analytic solution vector.  At the final time we compute the relative error
\begin{equation}
    \tilde{E} = \frac{||\mm{W}(T) - W^*(\cdot, T)||_\mm{H}}{||W^*(\cdot, T)||_\mm{H}},
\end{equation}
where $W^*(\cdot, T)$ is the analytic vector solution evaluated on the grid. Time stepping is done via an adaptive Runge-Kutta scheme with the time step sufficiently small so that the error is confined to that in space. Figures \ref{fig: conv1}-\ref{fig: conv4} illustrate  convergence for all cases at the theoretical rate.

\begin{figure*}[t!]
\centering
    \begin{subfigure}[t]{0.5\textwidth}
        \centering
        \includegraphics[height=1.7in]{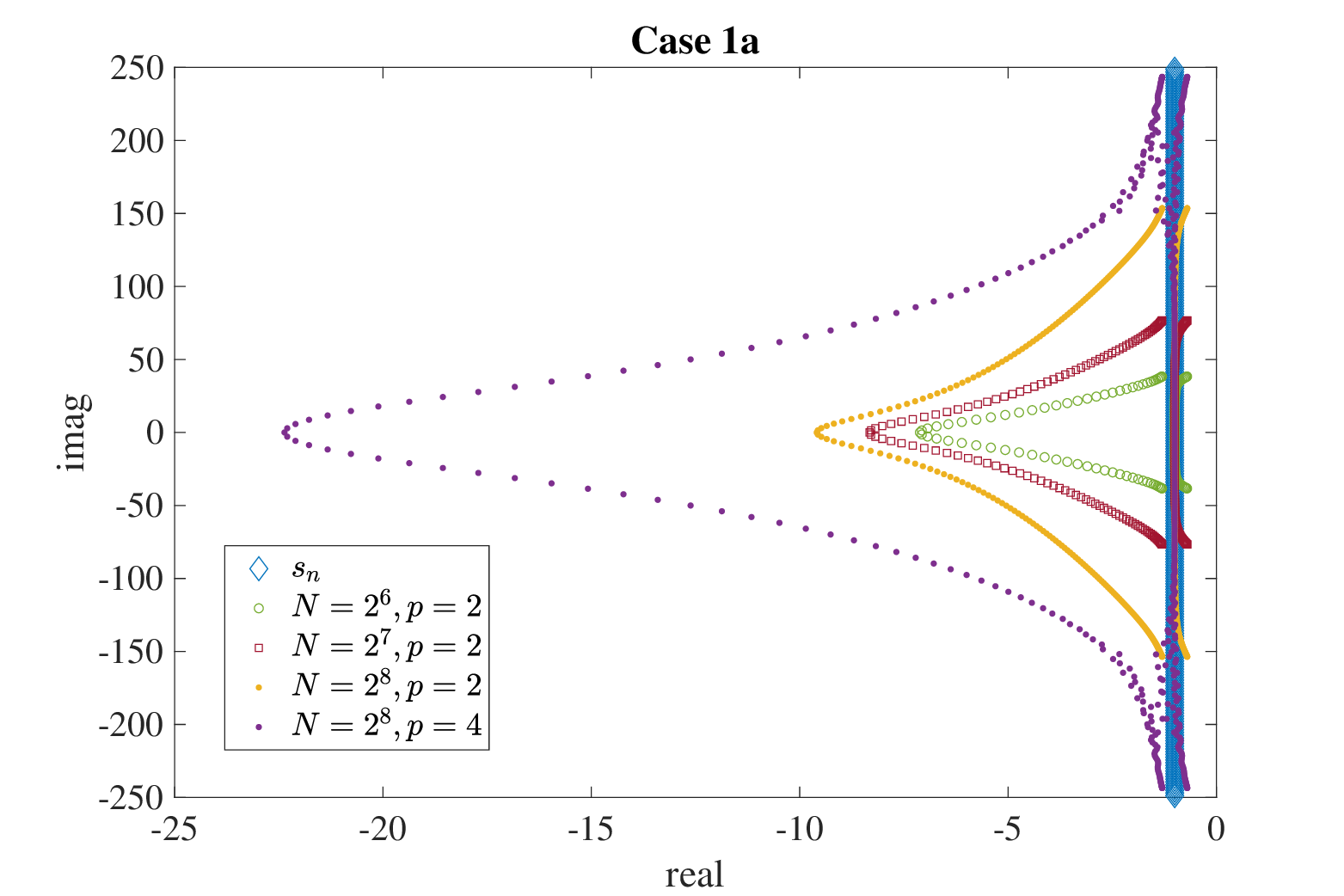}
        \caption{}
    \end{subfigure}%
    ~ \hspace{-5mm}
     \begin{subfigure}[t]{0.5\textwidth}
        \centering
        \includegraphics[height=1.7in]{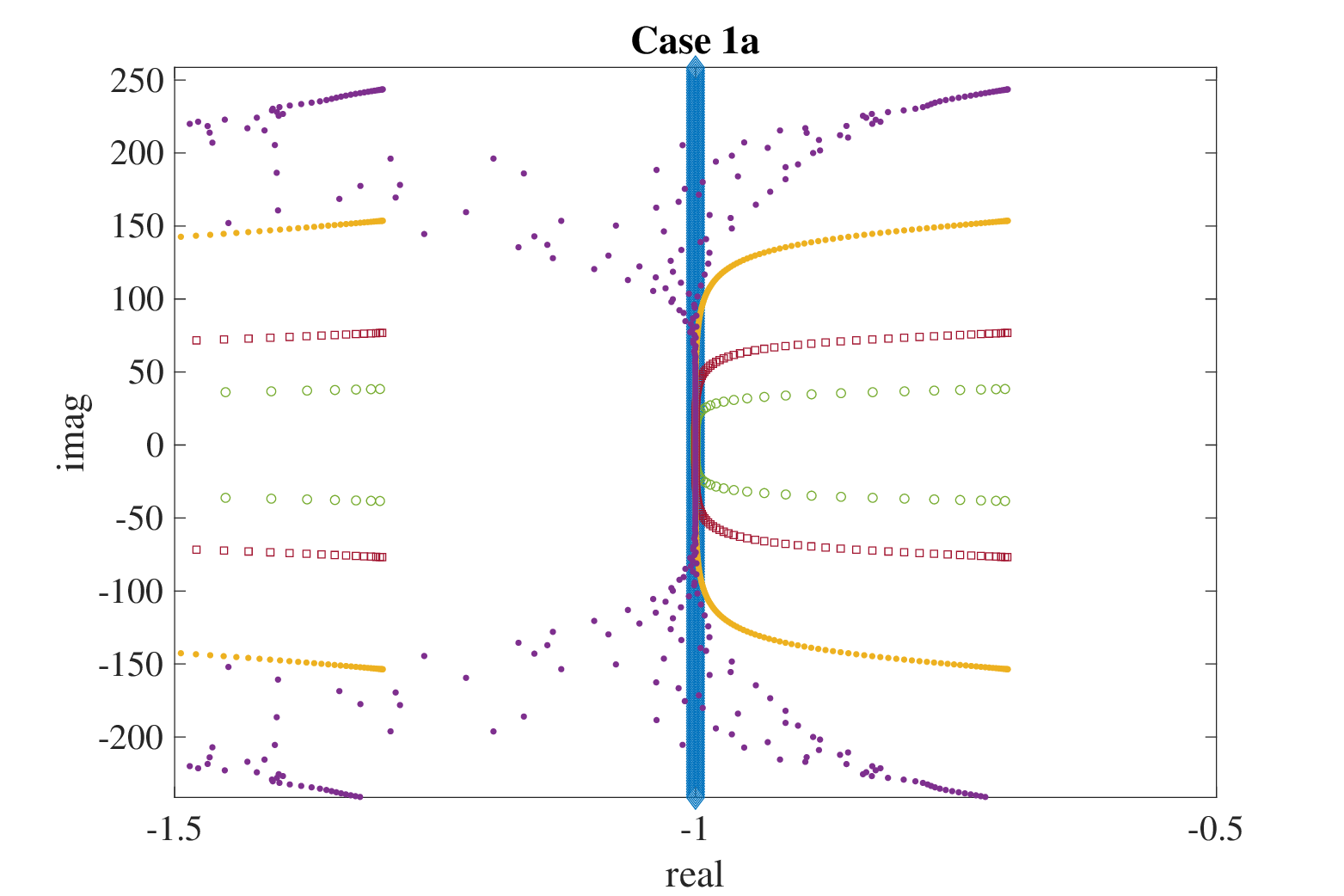}
        \caption{}
    \end{subfigure}
    \caption{(a) Analytic and numerical spectrum for Case 1a, where $a = -0.7, b = c = 0, d = -1.3$ with (b) zoom.}
    \label{fig: spec1a}
\end{figure*}

\begin{figure*}[t!]
\centering
    \begin{subfigure}[t]{0.5\textwidth}
        \centering
        \includegraphics[height=1.7in]{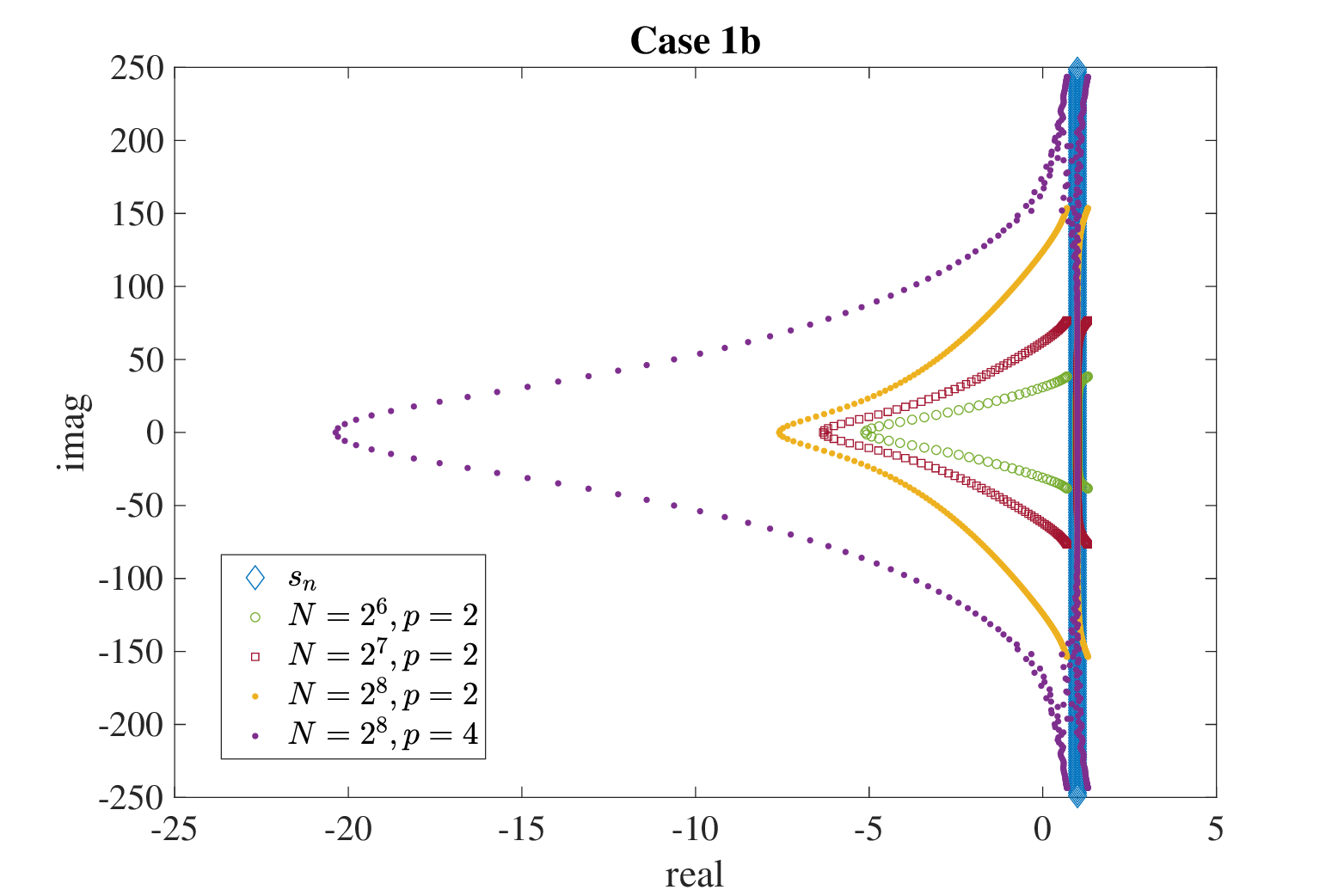}
        \caption{}
    \end{subfigure}%
    ~ \hspace{-5mm}
     \begin{subfigure}[t]{0.5\textwidth}
        \centering
        \includegraphics[height=1.7in]{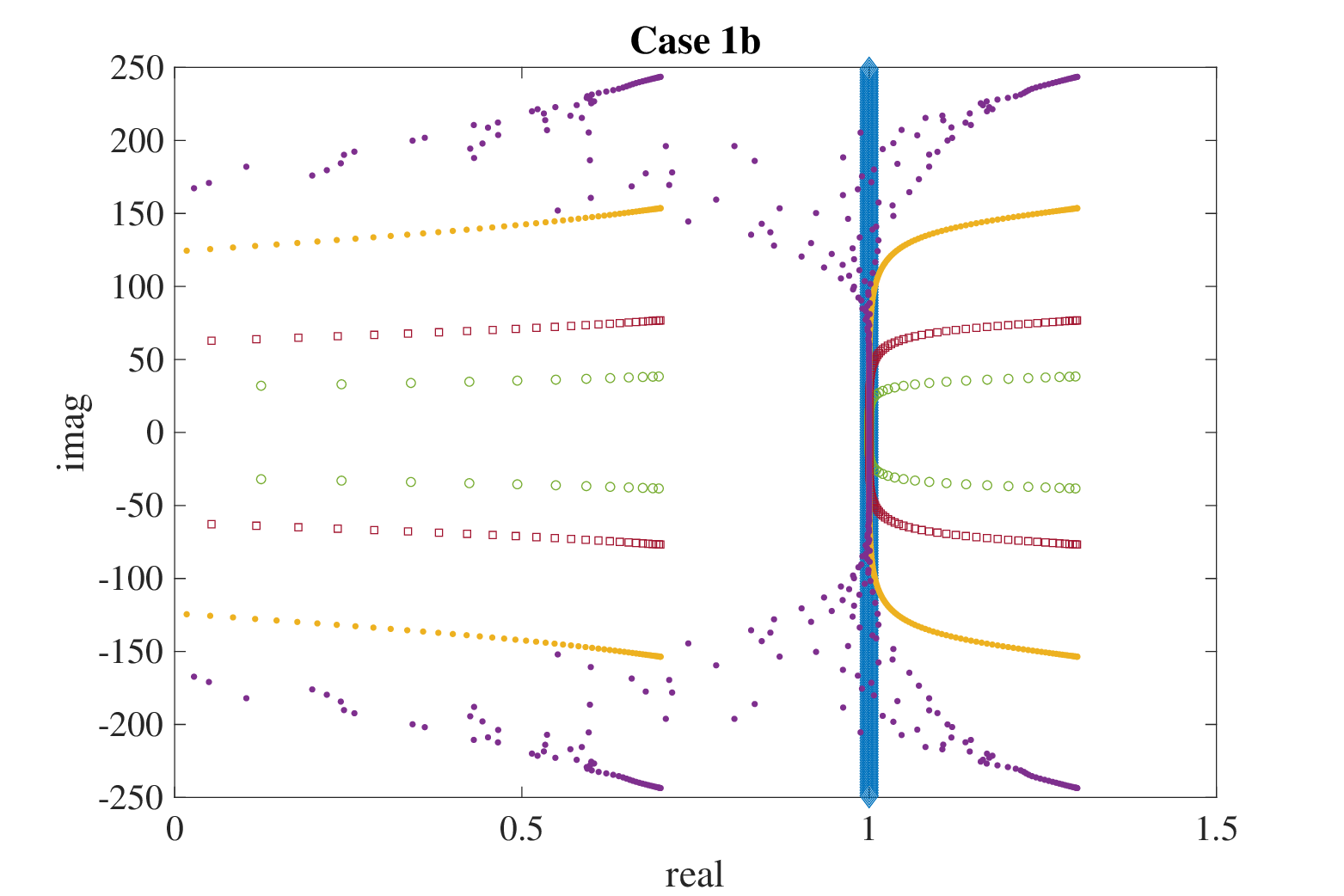}
        \caption{}
    \end{subfigure}
    \caption{(a) Analytic and numerical spectrum for case 1b, where $a = 0.7, b = c = 0, d = 1.3$ with (b) zoom.}
    \label{fig: spec1b}
\end{figure*}

\begin{figure*}[t!]
\centering
    \begin{subfigure}[t]{0.5\textwidth}
        \centering
        \includegraphics[height=1.7in]{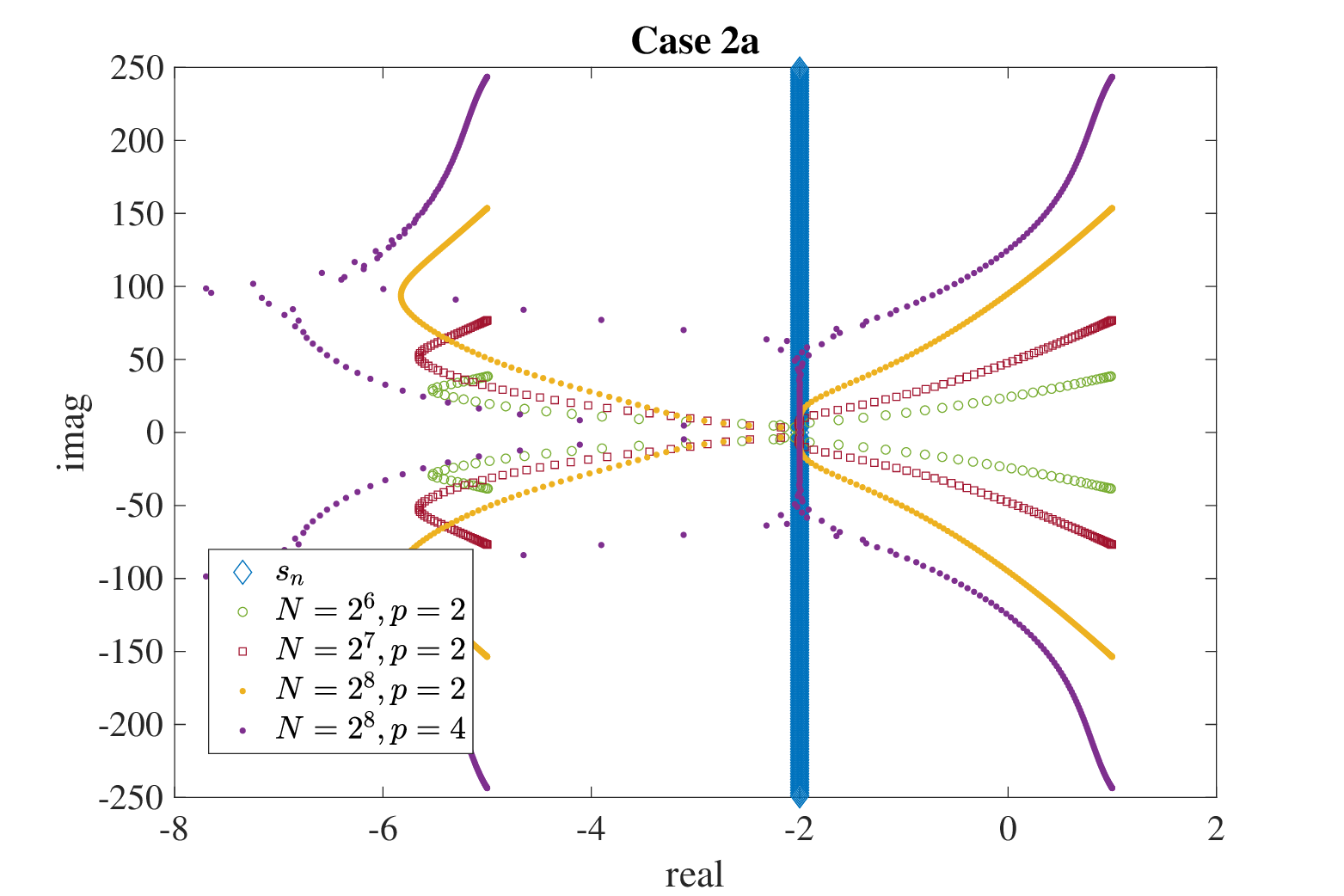}
        \caption{}
    \end{subfigure}%
    ~ \hspace{-5mm}
     \begin{subfigure}[t]{0.5\textwidth}
        \centering
        \includegraphics[height=1.7in]{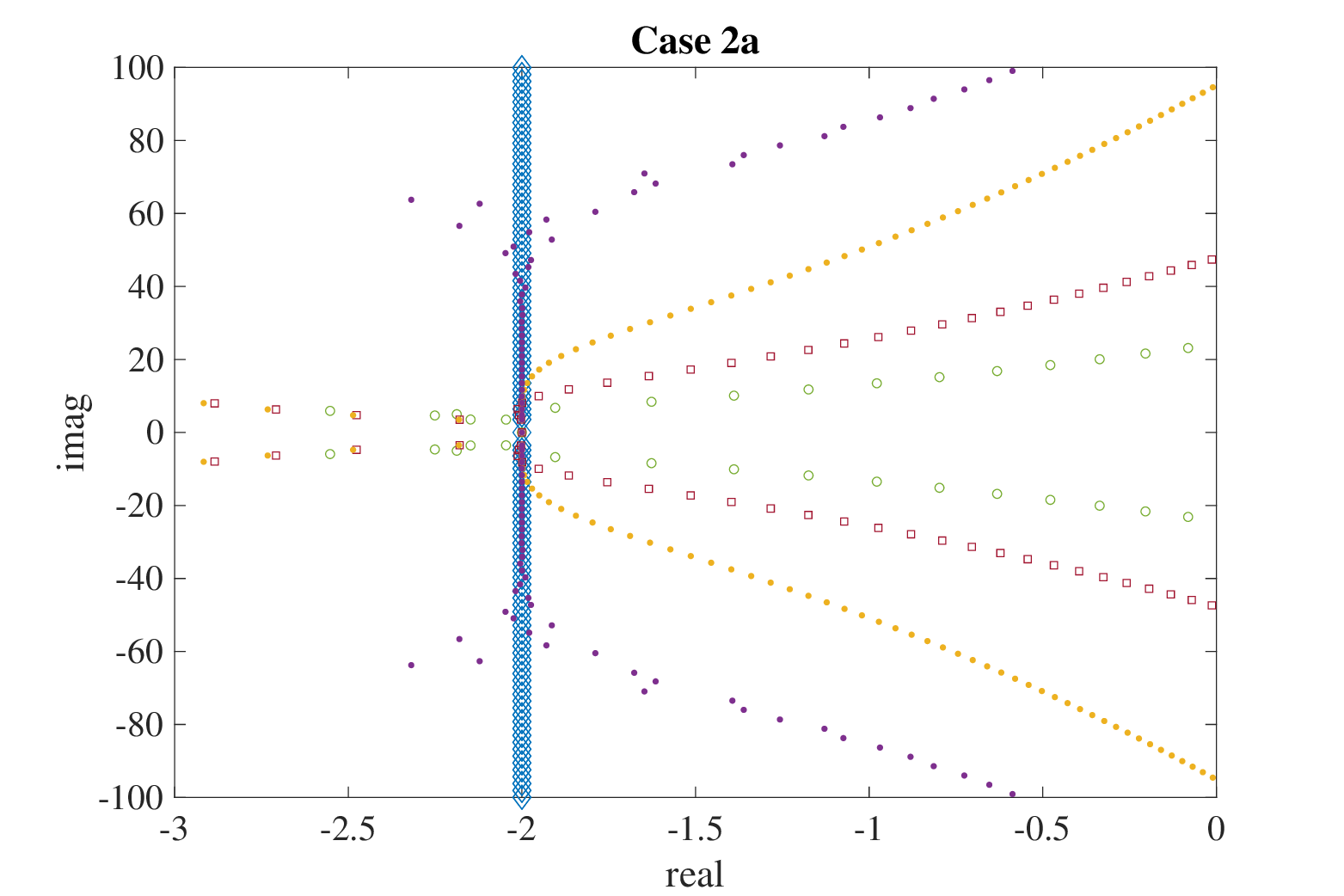}
        \caption{}
    \end{subfigure}
    \caption{(a) Analytic and numerical spectrum for case 2a, where $a = 1, b = -3, c = 3, d = -5$ with (b) zoom.}
    \label{fig: spec2a}
\end{figure*}

  \begin{figure*}[t!]
\centering
    \begin{subfigure}[t]{0.5\textwidth}
        \centering
        \includegraphics[height=1.7in]{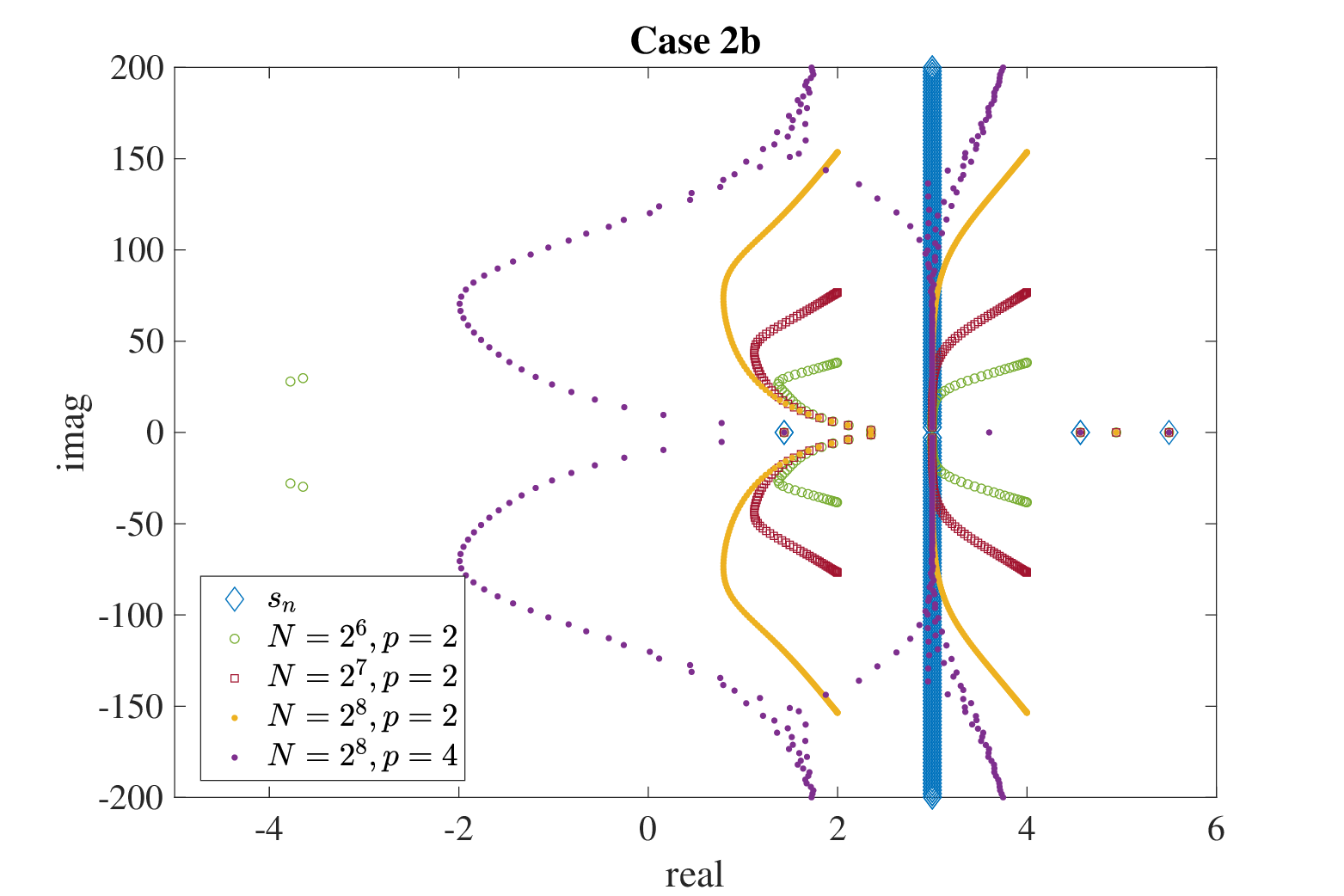}
        \caption{}
    \end{subfigure}%
    ~ \hspace{-5mm}
     \begin{subfigure}[t]{0.5\textwidth}
        \centering
        \includegraphics[height=1.7in]{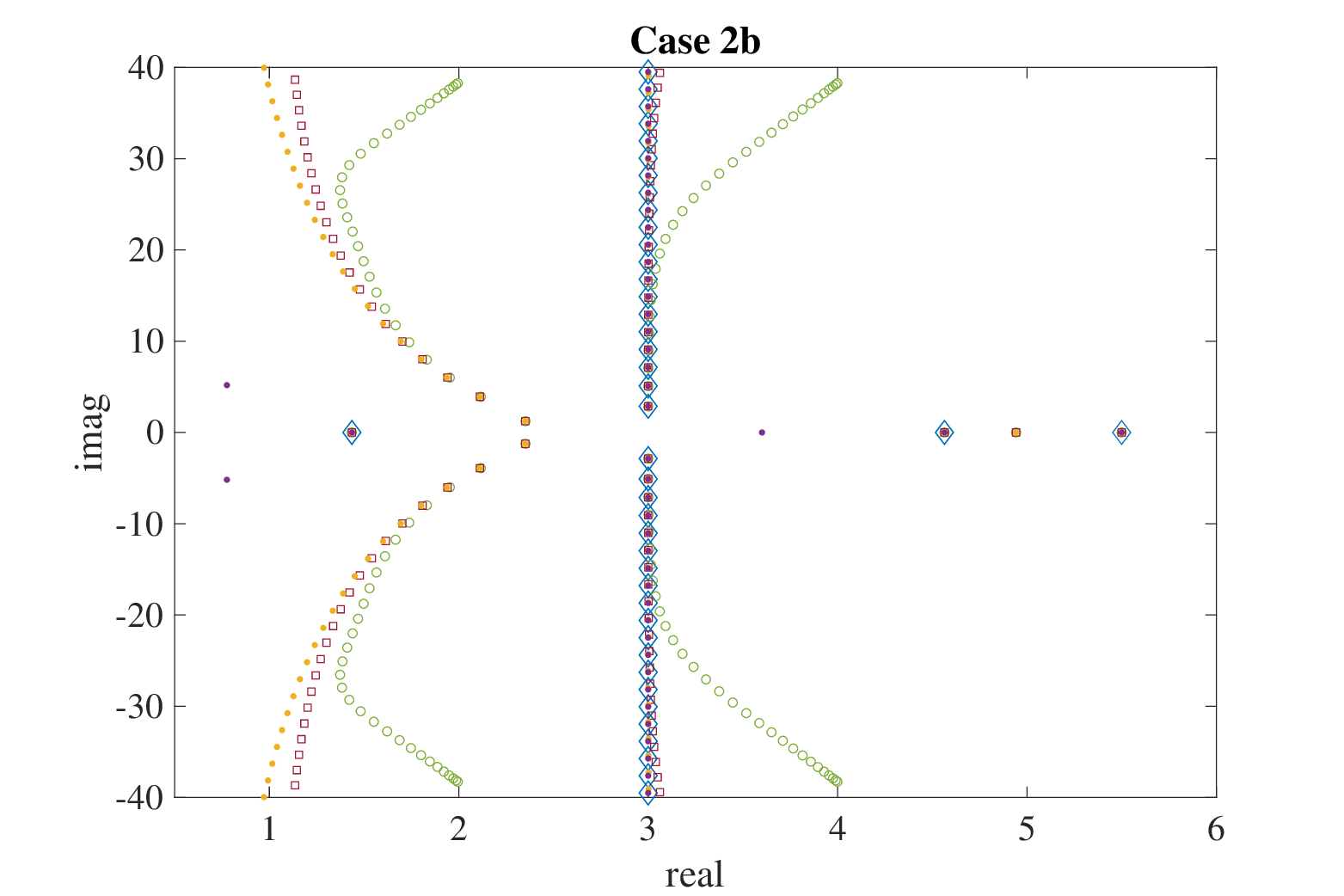}
        \caption{}
    \end{subfigure}
    \caption{(a) Analytic and numerical spectrum for case 2b, where $a = 2, b = 3, c = 2, d = 4$ with (b) zoom.}
    \label{fig: spec2b}
\end{figure*}

  \begin{figure*}[t!]
\centering
    \begin{subfigure}[t]{0.5\textwidth}
        \centering
        \includegraphics[height=1.7in]{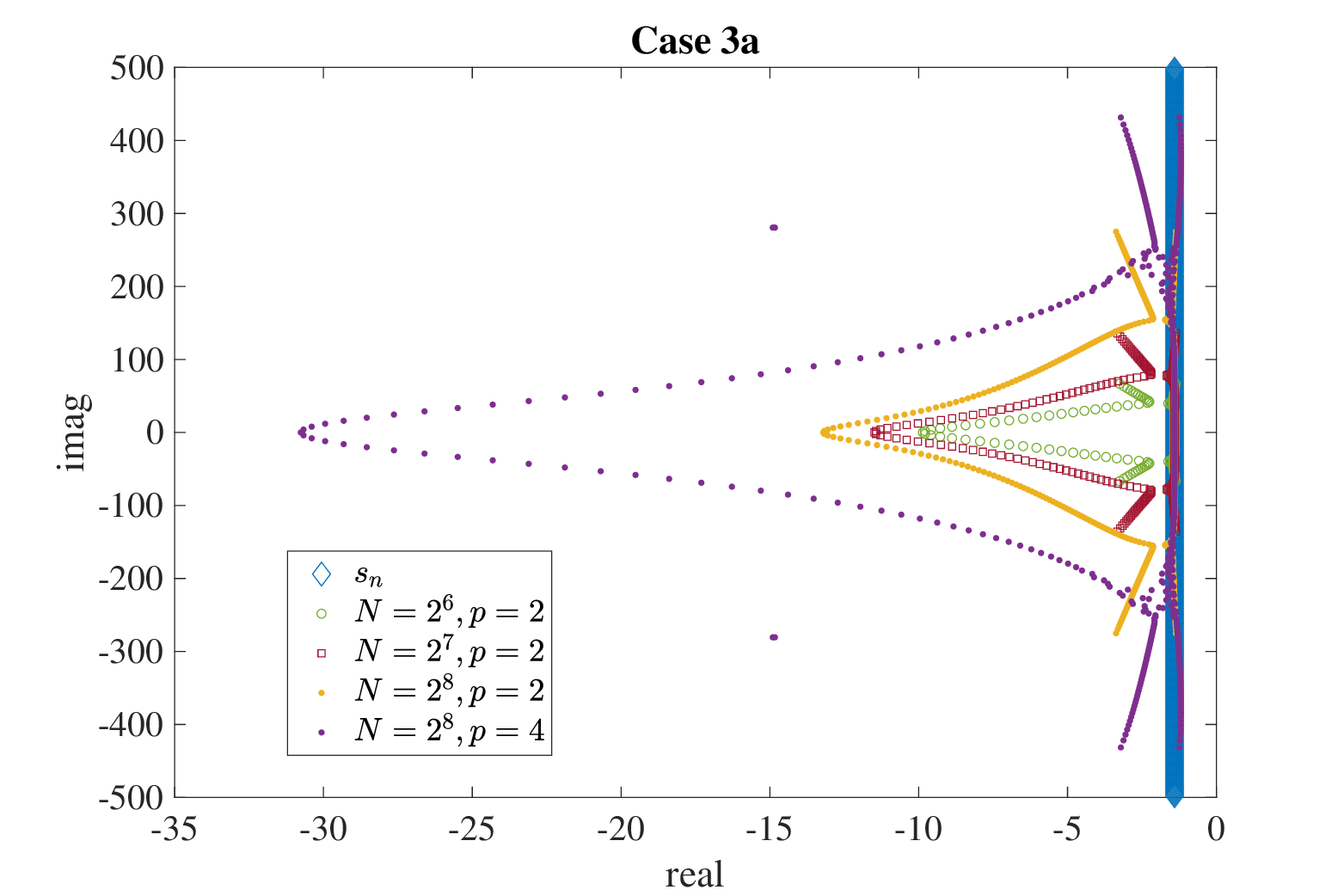}
        \caption{}
    \end{subfigure}%
    ~ \hspace{-5mm}
     \begin{subfigure}[t]{0.5\textwidth}
        \centering
        \includegraphics[height=1.7in]{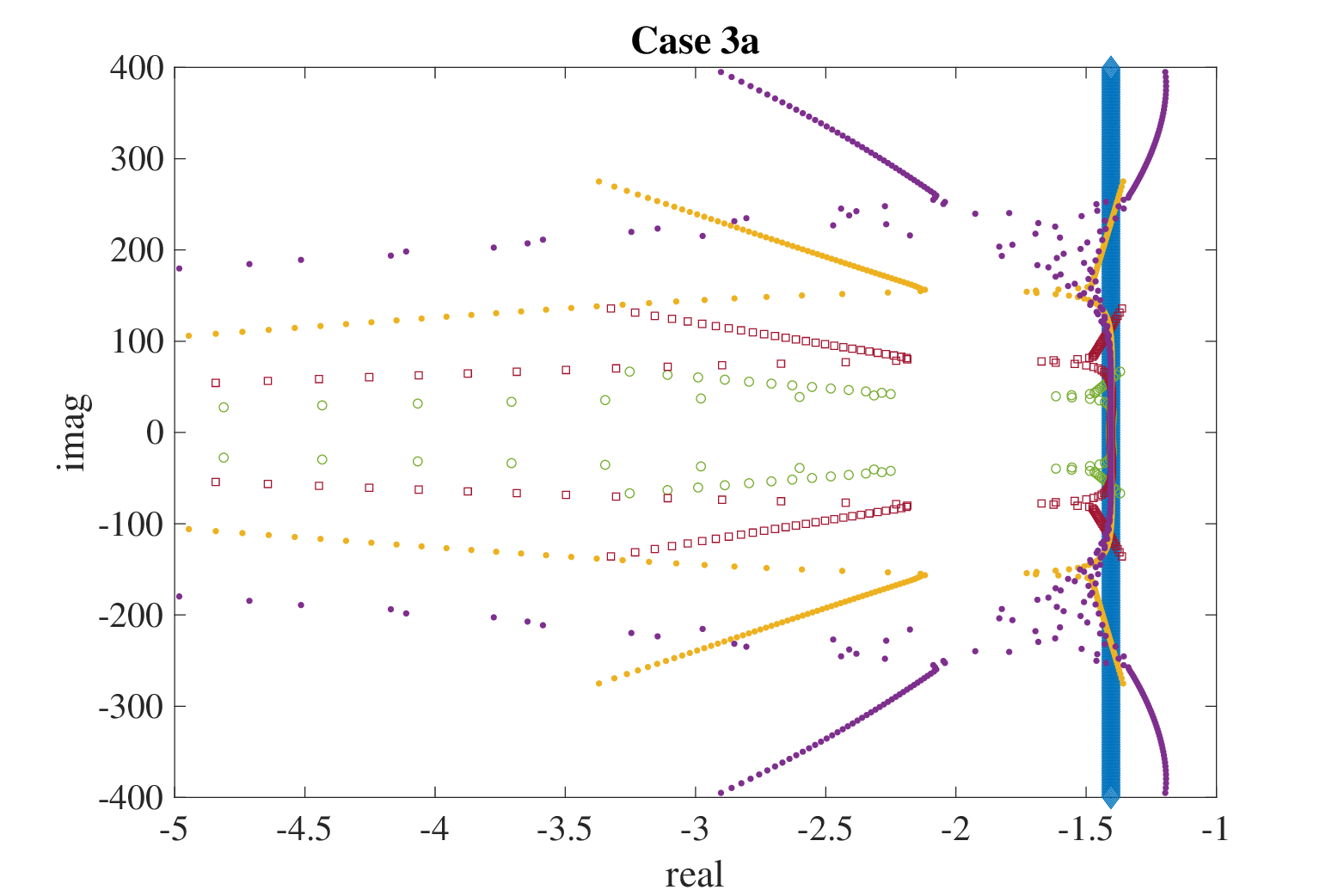}
        \caption{}
    \end{subfigure}
    \caption{(a) Analytic and numerical spectrum for case 3a, where $a(x) = -0.7 - x, d(x) = -1.3 + 0.1x$ with (b) zoom.}
    \label{fig: spec3a}
\end{figure*}

\begin{figure*}[t!]
\centering
    \begin{subfigure}[t]{0.5\textwidth}
        \centering
        \includegraphics[height=1.7in]{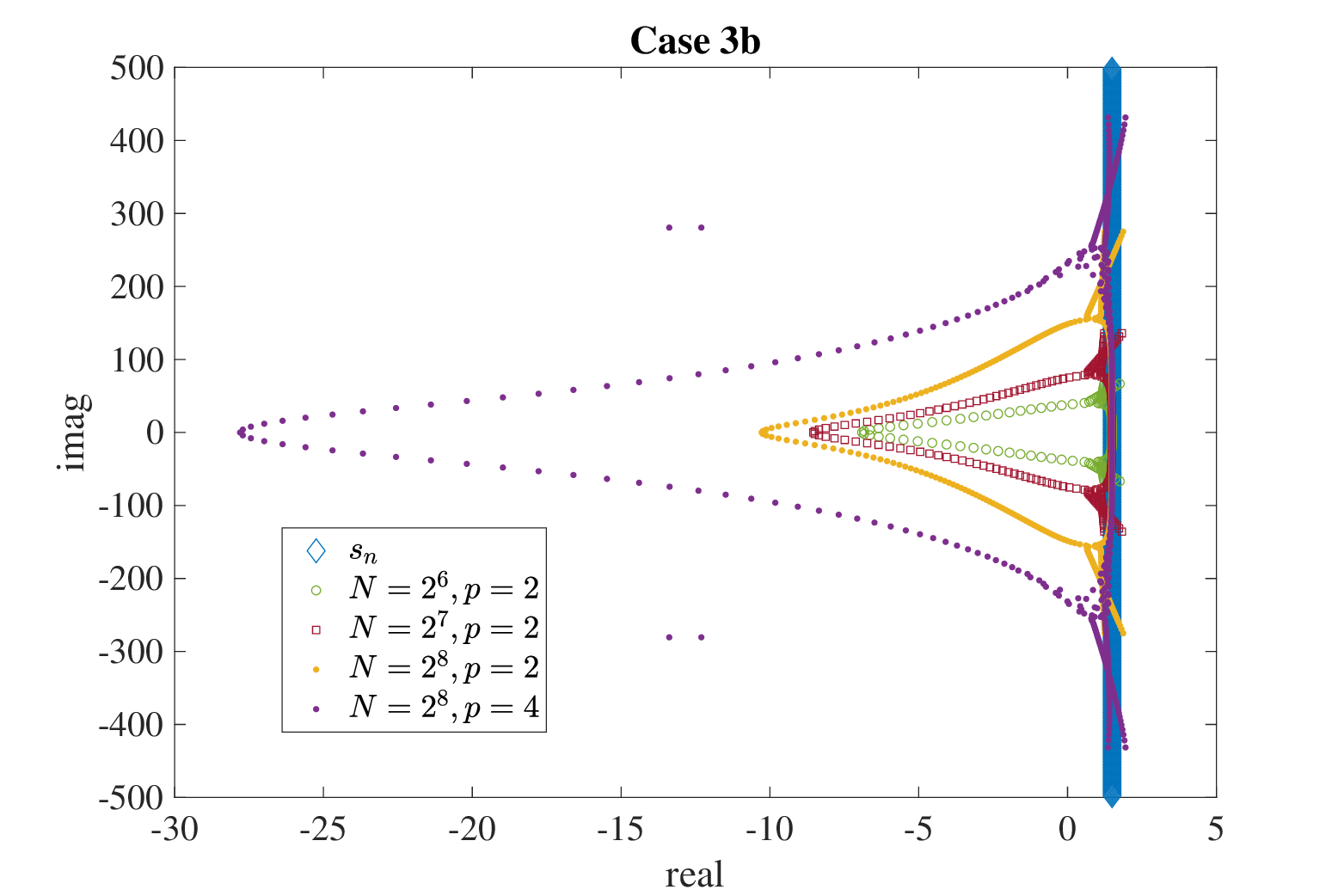}
        \caption{}
    \end{subfigure}%
    ~ \hspace{-5mm}
     \begin{subfigure}[t]{0.5\textwidth}
        \centering
        \includegraphics[height=1.7in]{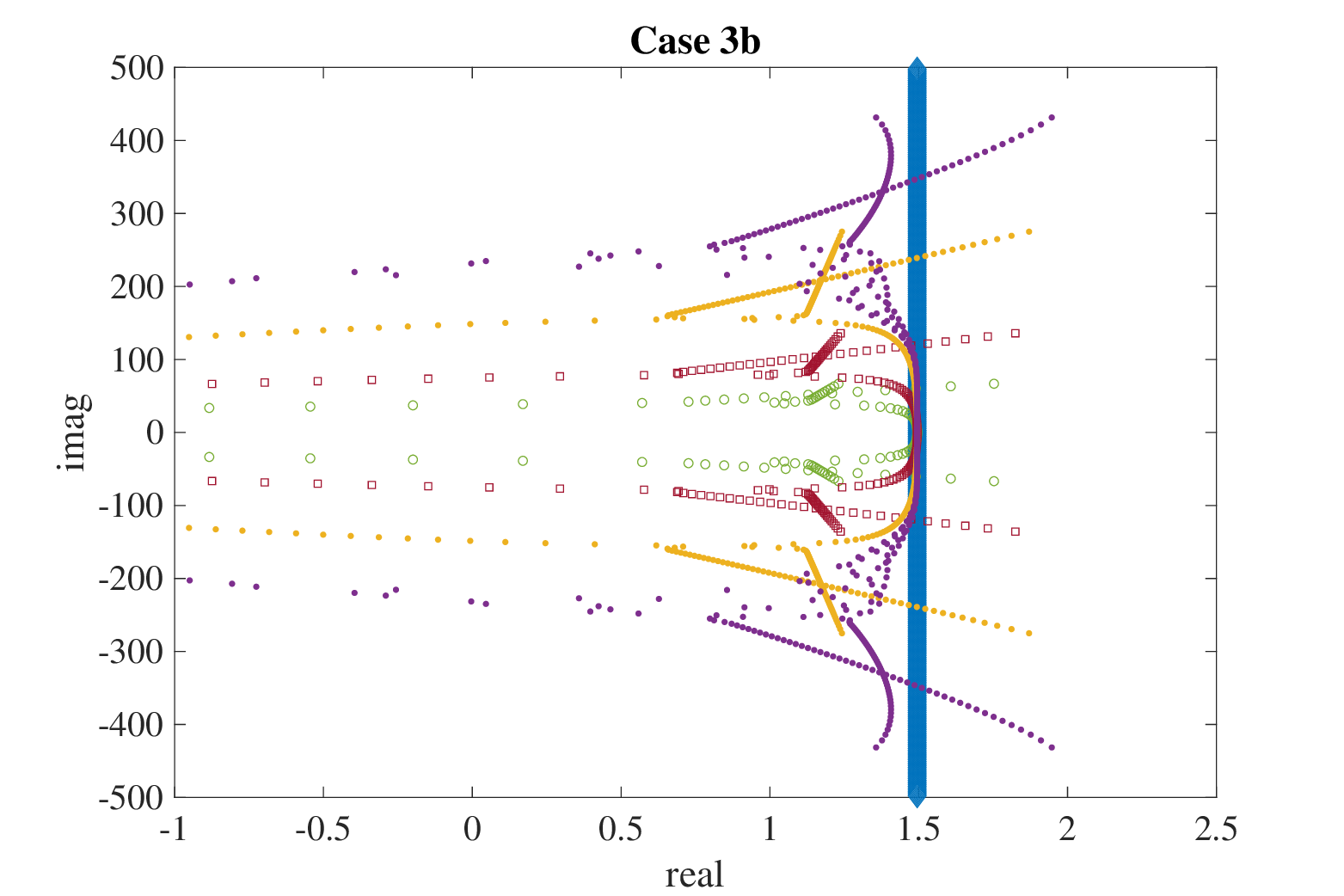}
        \caption{}
    \end{subfigure}
    \caption{(a) Analytic and numerical spectrum for case 3b, where $a(x) = 0.7 + x, d(x) = 1.3 + 0.1x$ with (b) zoom.}
    \label{fig: spec3b}
\end{figure*}

\begin{figure*}[t!]
\centering
    \begin{subfigure}[t]{0.5\textwidth}
        \centering
        \includegraphics[height=1.7in]{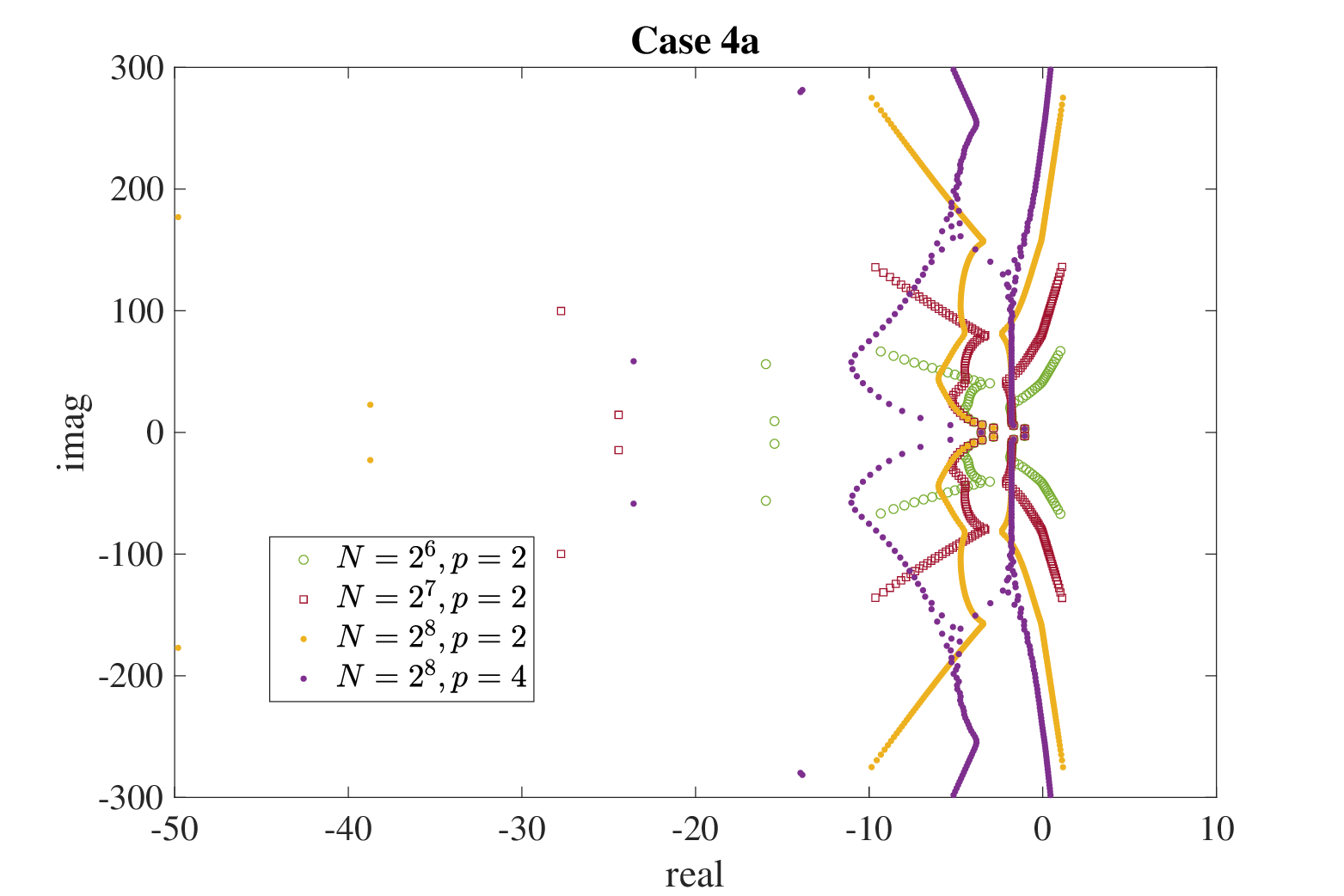}
        \caption{}
    \end{subfigure}%
    ~ \hspace{-5mm}
     \begin{subfigure}[t]{0.5\textwidth}
        \centering
        \includegraphics[height=1.7in]{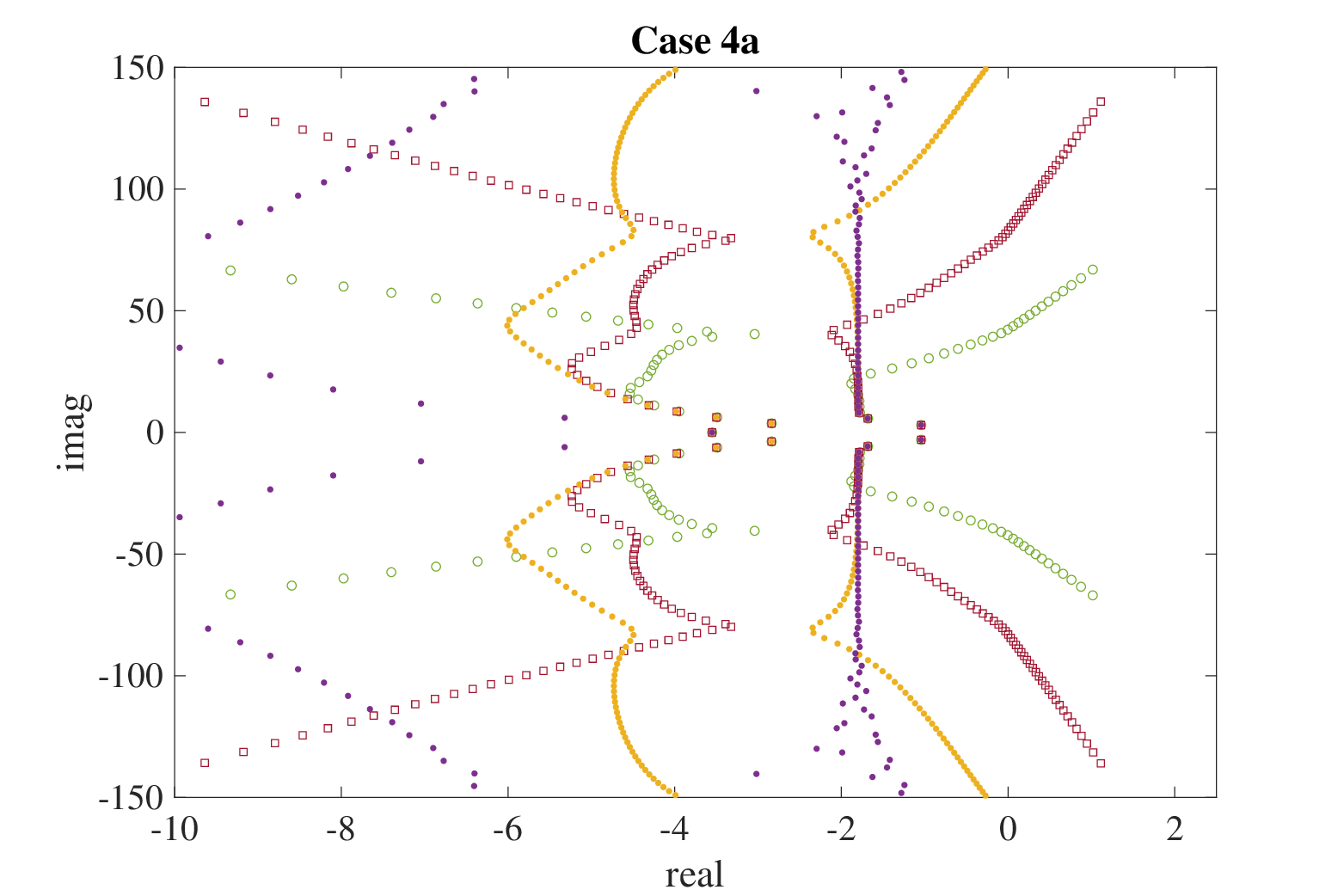}
        \caption{}
    \end{subfigure}
    \caption{(a) Analytic and numerical spectrum for case 4a, where $a(x) = x, b(x) = -3x, c(x) = 3x, d(x) = -5x$ with (b) zoom.}
    \label{fig: spec4a}
\end{figure*}

\begin{figure*}[t!]
\centering
    \begin{subfigure}[t]{0.5\textwidth}
        \centering
        \includegraphics[height=1.7in]{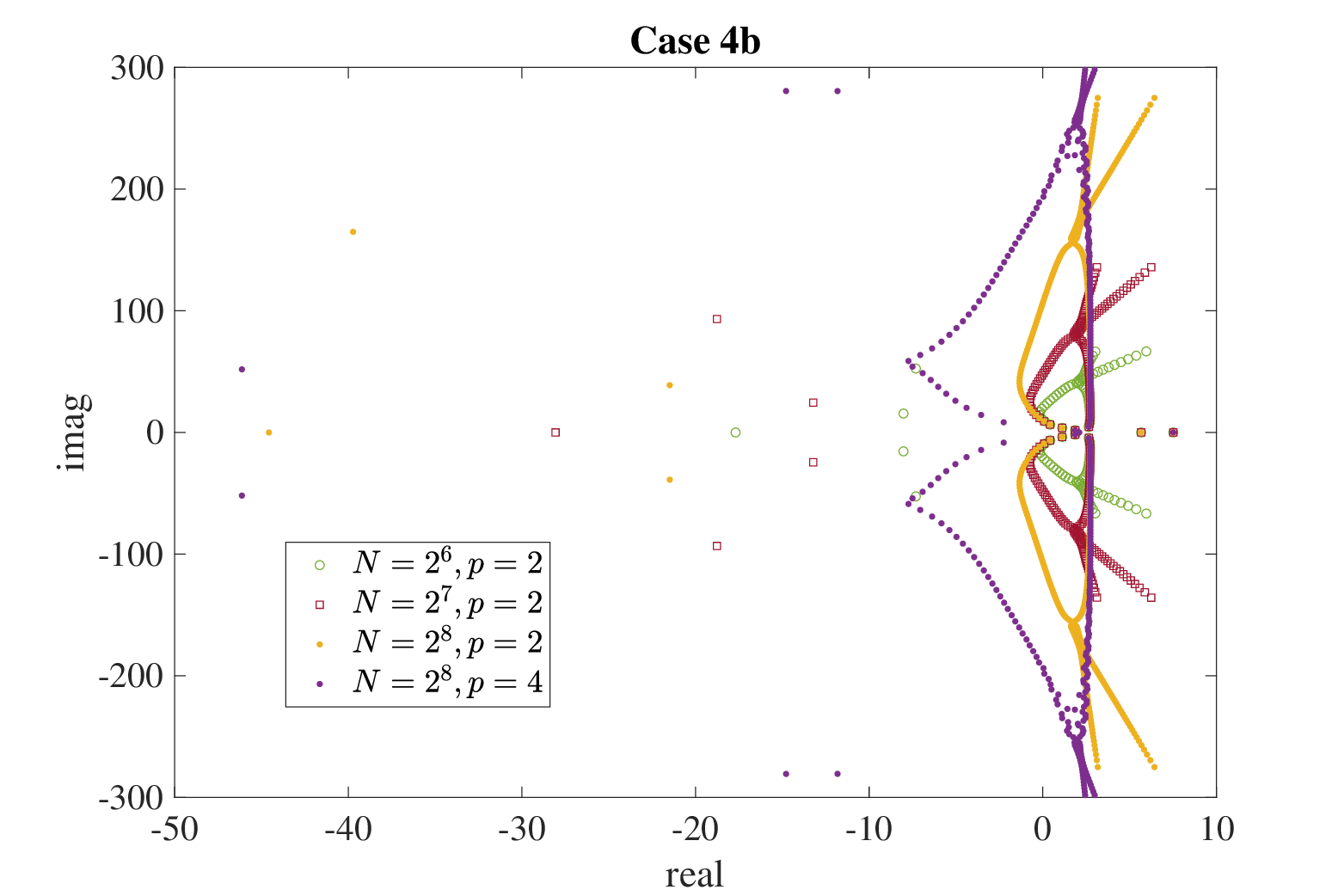}
        \caption{}
    \end{subfigure}%
    ~ \hspace{-5mm}
     \begin{subfigure}[t]{0.5\textwidth}
        \centering
        \includegraphics[height=1.7in]{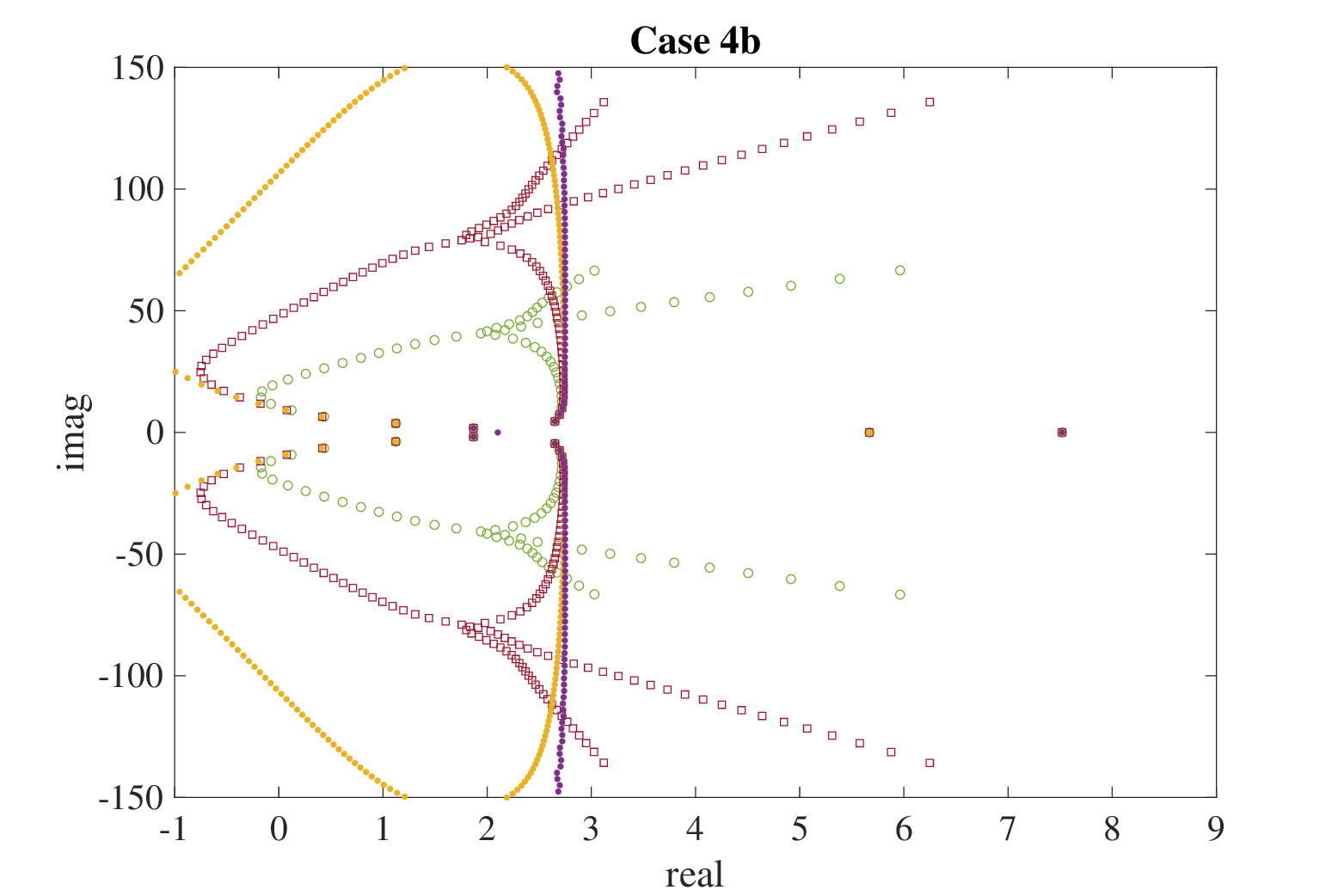}
        \caption{}
    \end{subfigure}
    \caption{(a) Analytic and numerical spectrum for case 4b, where $a(x) = 2 + x, b(x) = 3x, c(x) = 2x, d(x) = -1 + 4x$ with (b) zoom.}
    \label{fig: spec4b}
\end{figure*}

\section{Numerical Spectra Approximations}
\label{sec:discrete_spec}
Applying the Laplace transform to discretization \eqref{eqn: discrete2} yields the eigenvalue problem
\begin{equation}
    s\mm{W} = \mathcal{D}_h\mm{W},
\end{equation}
where $\mathcal{D}_h = (-\mm{\Lambda}_D  + \tilde{\mm{B}} + \mm{P})$, whose eigenvalues correspond to the numerical spectrum.  For Cases 1-4 we explore how well this numerical spectrum approximates the analytic spectrum, the latter of which was derived in section \ref{sec:cont_spec}. Ideally the numerical spectra does not overly predict the rate of temporal evolution in the physical problem, which would be the case if the numerical spectra converge to the analytic from the left, a feature that indicates strict stability of dissipative type of the numerical scheme \cite{Nord2006, Erickson2019}. This is important for long-time calculations so that high-frequency errors do not grow and destroy the accuracy. However we are primarily interested in determining when numerical spectra correctly identify positive modes of growth.

For Case 1a the analytic spectrum is $s_n = -1 + 0.6\pi i n, \quad n \in \mathbb{Z}$.   In  \Cref{fig: spec1a} we plot this along with the numerical spectrum for increasing number of grid points $N$, with $p = 2$ and also for one case of $p = 4$. The physical problem dictates exponential decay since $\text{Re}(s_n) < 0$, and the left figure shows that the numerical spectra lies entirely in the left half of the complex plane. However, the zoom on the right reveals that some of the numerical eigenvalues lay to the right of $s_n$, indicating that dissipative strict stability is not obtained, a feature that persists for all the cases we consider in this work. However, as the mesh is refined or if $p$ is increased, more numerical eigenvalues align with $s_n$, indicating convergence of the numerical spectrum.  We also observe that the real part of the numerical spectrum is bounded above by $a = -0.7$.  For Case 1b, the analytic spectrum is given by \begin{align}
  s_n = 1 + 0.6\pi i n, \quad n \in \mathbb{Z},  
\end{align} which correspond to physical growth. \Cref{fig: spec1b} reflects similar findings to that of Case 1a, with the real part of the numerical spectrum bounded above by $d = 1.3$. 

For Case 2a, we have the real root $s_0 = -2$ as well as the complex roots (with non-zero imaginary part)
\begin{align}
 s^\pm_n = -2 \pm \sqrt{-9 - 0.36\pi^2 n^2}, \quad n = 1, 2, ...,   
\end{align} with real part $\text{Re}(s^\pm_n) = -2$ corresponding to physical decay. \Cref{fig: spec2a} shows that some of the numerical eigenvalues have positive real part that remains bounded above by $a = 1$. This was a surprising finding for a problem with physical decay, but a stable method still allows for some exponential growth, as long as the solution is still bounded by the data of the problem \cite{bertil_2008}. Moreover, the numerical eigenvalues tend towards $s_n$ as the mesh is refined or $p$ increases.  To the left of $s_n$, the numerical spectra corresponding to $p = 2$ appears to persist with mesh refinement, which at first led us to believe that physical eigenvalues existed here. However, these do not persist with increasing $p$, rendering them numerical features.

For case 2b, we have the real root $s_0 = 5.5$ as well as
\begin{align*}
 s^\pm_n = 3 \pm \sqrt{6 - 0.36\pi^2 n^2}, \quad n = 1, 2, ...   
\end{align*} which generates real roots $s_1^+ \approx 4.5$, $s_1^- \approx 1.4$ with the remaining being complex roots with non-zero imaginary part and $\text{Re}(s_n) = 3$, indicating physical growth.  All three real eigenvalues are captured by the discretization, as illustrated in Figure \ref{fig: spec2b} and the discretization appears to converge (again from the right, although the real part of the eigenvalues containing non-zero imaginary parts are bounded above by $d = 4$). Note that in Cases 1 and 2, the real part of the numerical spectrum (for eigenvalues with non-zero imaginary part) is bounded above by the largest diagonal element of $\tilde{B}$; the proof that this is always the case is left for future work. 

For Case 3a, all the eigenvalues lay on a vertical line at $\text{Re}(s_n) \approx -1.4$, as evident in \Cref{fig: spec3a}. The numerical eigenvalues all have negative real part (indicating a stable discretization) and appear to converge to $s_n$ from the right.  $a(x) = -0.7 - x$ attains its maximum value of $-0.7$ on $x \in [0, 2]$ at $x = 0$ and $d(x) = -1.3 + 0.1x$ attains a maximum of $-1.1$ at $x = 2$, the latter of which appears to provide an upper bound on the real part of the numerical spectrum.  The observation that numerical eigenvalues with non-zero imaginary part are bounded above by the maximum value attained on the diagonal of $\tilde{B}$ persists in the remaining cases with variable coefficients. For Case 3b, $\text{Re}(s_n) \approx 1.49$, with some numerical eigenvalues to the right, which appear to have real part bounded above the maximum value that $a(x)$ attains on $x \in [0, 2]$, although the bound is not as tight as in the constant coefficient case, see \Cref{fig: spec3b}. For Case 4 we do not have an analytic spectrum with which to compare. However, the numerical spectrum shown in \Cref{fig: spec4a} suggests negative growth, as eigenvalues that lay to the right of zero do not persist with mesh refinement or increasing order of accuracy.  For Case 4b, the spectrum shown in \Cref{fig: spec4a} suggests a physical instability is present as numerical eigenvalues with positive real part persist with mesh refinement and increasing $p$, namely at the complex numbers with $\text{Re}(z) \approx 2, 2.7$ and $7.5$. \section{Conclusions}

\label{sec:conclusions}

In this work we derived analytic spectra and solutions to systems of hyperbolic partial differential equations with variable coefficients that arise in perturbed problems, for example when linearizing a nonlinear system around a steady-state solution.  We developed a provably stable semi-discretization based on high-order SBP-SAT finite difference schemes in order to explore numerical investigations of physical instabilities.  While we did not obtain a strictly stable method of dissipative type (so that the growth rate of the discretization is bounded by that of the analytic) we did observe convergence of the numerical spectrum. 

We conclude by stating that caution should be taken when using discrete methods to approximate analytic spectra. A method that is provably stable is desirable, and demonstration of numerical convergence will ensure that positive eigenvalues are physical and not a result of a numerical instability. When the analytic spectra and/or solutions cannot be computed in closed form, discrete approximations to the spectra should exhibit convergence with both mesh-refinement and with higher-order spatial accuracy in order to gain assurance that one is well-capturing the underlying physics.

%% Refer following link for more details.
%% https://en.wikibooks.org/wiki/LaTeX/Mathematics
%% https://en.wikibooks.org/wiki/LaTeX/Advanced_Mathematics

%% Use a table environment to create tables.
%% Refer following link for more details.
%% https://en.wikibooks.org/wiki/LaTeX/Tables
%\begin{table}[t]%% placement specifier
%% Use tabular environment to tag the tabular data.
%% https://en.wikibooks.org/wiki/LaTeX/Tables#The_tabular_environment
%\centering%% For centre alignment of tabular.
%\begin{tabular}{l c r}%% Table column specifiers
%% Tabular cells are separated by &
%  1 & 2 & 3 \\ %% A tabular row ends with \\
%  4 & 5 & 6 \\
%  7 & 8 & 9 \\
%\end{tabular}
%% Use \caption command for table caption and label.
%\caption{Table Caption}\label{fig1}
%\end{table}

%% Use figure environment to create figures
%% Refer following link for more details.
%% https://en.wikibooks.org/wiki/LaTeX/Floats,_Figures_and_Captions
%\begin{figure}[t]%% placement specifier
%% Use \includegraphics command to insert graphic files. Place graphics files in 
%% working directory.
%\centering%% For centre alignment of image.
%\includegraphics{example-image-a}
%% Use \caption command for figure caption and label.
%\caption{Figure Caption}\label{fig1}
%% https://en.wikibooks.org/wiki/LaTeX/Importing_Graphics#Importing_external_graphics
%\end{figure}

\section*{Acknowledgements} B.A.E. was funded by National Science Foundation awards \#2036980  and \#2053372.

%%       \cite{<label>} ==> [1]

%%

%  \bibliographystyle{elsarticle-num} 
%  \bibliography{references.bib}

\end{document}